\documentclass[journal]{IEEEtran}
\usepackage{amssymb,amstext,amsmath,amsthm}
\usepackage{multicol}
\usepackage{epsfig,float,lmodern}
\usepackage[compress]{cite}

\newtheorem{example}{Example}
\newtheorem{remark}{Remark}
\newtheorem{definition}{Definition}

\begin{document}

\title{
 Hold-in, pull-in, and lock-in ranges of PLL circuits:
 rigorous mathematical definitions \\ and limitations of classical theory.
}

\author{
  Leonov~G.A., Kuznetsov~N.V., Yuldashev~M.V., Yuldashev~R.V.
\thanks{
}
\thanks{
Saint-Petersburg State University, Russia;
University of Jyv\"{a}skyl\"{a}, Finland
(corr. author e-mail:nkuznetsov239@gmail.com).
}
\thanks{
  IEEE Transactions on Circuits and Systems--I: Regular Papers;
  DOI 10.1109/TCSI.2015.2476295 (accepted)
}
}

\maketitle

\begin{abstract}
The terms hold-in, pull-in (capture), and lock-in ranges are widely used
by engineers for the concepts of frequency deviation ranges
within which PLL-based circuits can achieve lock under various additional conditions.
Usually only non-strict definitions are given for these concepts in engineering literature.
After many years of their usage,
F.~Gardner in the 2nd edition of his well-known work,
{\emph{Phaselock Techniques}},
wrote
{\it{``There is no natural way to define exactly any unique lock-in frequency}''}
and {\it{``despite its vague reality, lock-in range is a useful concept}''}.
Recently these observations have led to the following advice given in
a handbook on synchronization and communications
``{\it{We recommend that you check these definitions carefully
before using them}}''\cite[p.49]{KiharaOE-2002}.
In this survey
an attempt is made to discuss and fill some of the gaps identified
between mathematical control theory,  the theory of dynamical systems
and the engineering practice of phase-locked loops.
It is shown that, from a mathematical point of view,
in some cases the hold-in and pull-in ``ranges''
may not be the intervals of values but a union of intervals and
thus their widely used definitions require clarification.
Rigorous mathematical definitions for the hold-in, pull-in, and lock-in ranges
are given.
An effective solution for the problem on the unique definition of the lock-in frequency,
posed by Gardner, is suggested.
\end{abstract}


\section{Introduction}
\IEEEPARstart{T}{he} phase-locked loop based circuits (PLL)
are widely used in various applications. A PLL is essentially a nonlinear
control system and its nonlinear analysis is a challenging task.
Much engineering writing is devoted to the study
of PLL-based circuits and the various characteristics for their stability
(see, e.g. a rather comprehensive bibliography of pioneering works in \cite{LindseyT-1973}).
An important engineering characteristic of PLL
is a set of parameters' values for which a PLL achieves lock.
In the classical books on PLLs \cite{Gardner-1966,Viterbi-1966,ShahgildyanL-1966},
published in 1966, such concepts as hold-in, pull-in, lock-in,
and other frequency ranges for which PLL can achieve lock,
were introduced. They are widely used nowadays
(see, e.g. contemporary engineering literature
\cite{Gardner-2005-book,Best-2007,Kroupa-2012} and other publications).
Usually in engineering literature only non-strict definitions are given
for these concepts.
F.~Gardner in 1979\footnote{A year later, in 1980, F.Gardner was elected IEEE Fellow
``for contributions to the understanding and applications of phase lock loops''.}
in the 2nd edition of his well-known work, \emph{Phaselock Techniques},
formulated the following problem \cite[p.70]{Gardner-1979-book}
(see also the 3rd edition \cite[p.187-188]{Gardner-2005-book}):
 ``{\it{There is no natural way to define exactly any unique lock-in frequency}}''.
The lack of rigorous explanations led
to the paradox:
 ``{\it{despite its vague reality, lock-in range is a useful concept}}''
\cite[p.70]{Gardner-1979-book}.
Many years of using definitions based on the above concepts
has led to the advice given in a handbook on synchronization and
communications, namely to check the definitions carefully
before using them \cite[p.49]{KiharaOE-2002}.
\smallskip

In this paper it is shown that, from a mathematical point of view,
in some cases the hold-in and pull-in ``ranges''
may be not intervals of values but a union of intervals,
and thus their widely used definitions require clarification.
Next, rigorous mathematical definitions for the hold-in, pull-in, and lock-in ranges are given.
In addition we suggest an effective solution
for the problem of the unique definition of the lock-in frequency,
posed by Gardner.

\section{Classical nonlinear mathematical models
of PLL-based circuits in a signal's phase space} % a - inserted by the editor

\begin{figure}[H]
 \centering
 \includegraphics[width=0.73\linewidth]{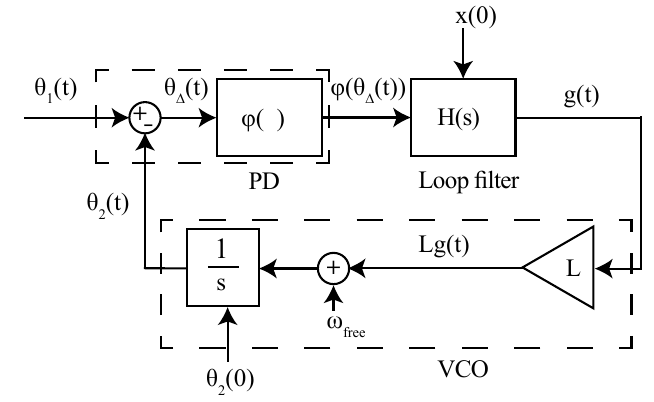}
 \caption{
   PLL-based circuit in a signal's phase space.
 }
 \label{costas-pll-sim-model}
\end{figure}

In classical engineering publications various analog PLL-based circuits
are represented in a \emph{signal's phase space}
(also named \emph{frequency-domain} \cite[p.338]{Davis-2011})
by the block diagram shown in Fig.~\ref{costas-pll-sim-model}.
Considering the corresponding mathematical model:
the Phase Detector (PD) is a nonlinear block;
the phases $\theta_{1,2}(t)$ of
the input (reference) and VCO signals
are PD block inputs and the output is
a function $\varphi(\theta_\Delta(t)) = \varphi(\theta_1(t)-\theta_2(t))$
named a \emph{phase detector characteristic},
where
\begin{equation}
  \label{theta_delta_def}
  \begin{aligned}
    & \theta_\Delta(t) = \theta_1(t) - \theta_2(t),
  \end{aligned}
\end{equation}
named the phase error.
The relationship between the input $\varphi(\theta_\Delta(t))$
and the output $g(t)$ of the linear filter (Loop filter) is as follows:
\begin{equation}\label{loop-filter}
 \begin{aligned}
 & \dot x = A x + b \varphi(\theta_\Delta(t)),
 \ g(t) = c^*x + h\varphi(\theta_\Delta(t)),
 \end{aligned}
\end{equation}
where $A$ is a constant matrix,
$x(t) \in \mathbb{R}^n$ the filter state,
$x(0)$ the initial state of filter,
$b$ and $c$ constant vectors, and h a number.
The filter transfer function has the form:\footnote{
  In the control theory the transfer function is often defined with the opposite sign
  (see, e.g. \cite{LeonovK-2014}): $H(s) = c^*(A-sI)^{-1}b-h.$
}
\begin{equation}
  H(s) = -c^*(A-sI)^{-1}b+h.
\end{equation}
A lead-lag filter \cite{Best-2007} (usually $H(0)=-c^*A^{-1}b+h=1$,
but $H(0)$ can also be any nonzero value
when an active lead-lag filter is used), or a PI filter ($H(0)$ is infinite)
is usually used as the filter.
The solution of \eqref{loop-filter} with initial data $x(0)$
(the filter output for the initial state $x(0)$) is as follows:
\begin{equation}\label{loop-filter-int}
 \begin{array}{c}
 g(t,x_0) = \alpha_0(t,x(0)) +
 \int\limits_0^t
 \gamma(t - \tau)\varphi(\theta_\Delta(\tau))
 {\rm d}\tau
 + h \varphi(\theta_\Delta(t)),
 \end{array}
\end{equation}
where $\gamma(t - \tau)=c^*e^{A(t-\tau)}b$ is the impulse response function of the filter
and $\alpha_0(t,x(0))= c^*e^{At}x(0)$ the zero input response
(natural response, i.e. when the input of the filter is zero).
The control signal $g(t)$ adjusts the VCO frequency to
the frequency of the input signal:
\begin{equation} \label{vco first}
   \dot\theta_2(t) = \omega_2(t) = \omega_2^{\text{free}} + Lg(t),
\end{equation}
where $\omega_2^{\text{free}}$ is the VCO free-running frequency
(i.e. for $g(t)\equiv 0$) and $L$ the VCO gain.
Nonlinear VCO models can be similarly considered, see, e.g. \cite{Margaris-2004,Suarez-2009}.
The frequency of the input signal (reference frequency) is usually assumed
to be constant:
\begin{equation}\label{omega1-const}
  \dot\theta_1(t) = \omega_1(t) \equiv \omega_1.
\end{equation}
The difference between the reference frequency and the VCO free-running frequency
is denoted as $\omega_\Delta^{\text{free}}$:
\begin{equation}
  \label{omega_delta_def}
  \begin{aligned}
    & \omega_\Delta^{\text{free}} \equiv \omega_1 - \omega_2^{\text{free}}.
  \end{aligned}
\end{equation}

By combining equations \eqref{theta_delta_def}, \eqref{loop-filter}, and \eqref{vco first}--\eqref{omega_delta_def}
a \emph{nonlinear mathematical model in the signal's phase space} is obtained
(i.e. in the state space: the filter's state $x$ and the difference between the signal's phases $\theta_\Delta$):
\begin{equation}\label{final_system}
 \begin{aligned}
   & \dot{x} = A x + b \varphi(\theta_{\Delta}), \\
   & \dot\theta_{\Delta} = \omega_{\Delta}^{\text{free}}
   - Lc^*x - Lh\varphi(\theta_{\Delta}).
 \end{aligned}
\end{equation}
Nowadays nonlinear model \eqref{final_system} is widely used
(see, e.g. \cite{Abramovitch-2002,Abramovitch-2004,Best-2007})
to study acquisition processes of various circuits.
The model can be obtained from the corresponding model in \emph{the signal space}
(called also \emph{time-domain} \cite[p.329]{Davis-2011})
by averaging under certain conditions \cite{KrylovB-1947,KudrewiczW-2007,LeonovKYY-2012-TCASII,LeonovK-2014,LeonovKYY-2015-SIGPRO},
a rigorous consideration of which is often omitted
(see, e.g. classical books \cite[p.12,15-17]{Viterbi-1966}, \cite[p.7]{Gardner-1966})
while their violation may lead to unreliable results
(see, e.g. \cite{KuznetsovKLNYY-2015-ISCAS,BestKKLYY-2015-ACC}).

Usually the PD characteristic is an odd function
(e.g. a PD realization such as a multiplier,
JK-flipflop, EXOR, PFD,
and other elements \cite{Best-2007}).
Note that the PD characteristic $\varphi(\theta_{\Delta})$ depends
on the waveforms of the considered
signals \cite{LeonovKYY-2012-TCASII,LeonovKYY-2015-SIGPRO}).
For the classical PLL with sinusoidal signals
and a two-phase PLL we have $\varphi(\theta_\Delta)=\frac{1}{2}\sin(\theta_\Delta)$,
for the classical BPSK Costas loop with ideal low-pass filters and a two-phase Costas loop
we have $\varphi(\theta_\Delta)=\frac{1}{8}\sin(2\theta_\Delta)$.

Classical PD characteristics are bounded piecewise smooth $2\pi$ periodic
functions\footnote{
If $\varphi(\theta_\Delta(t))$ has another period
(e.g. $\pi$ for the Costas loop models),
it has to be considered in the further discussion instead of $2\pi$.
}:
\[
  \varphi(\theta_\Delta(t)+2\pi k)=\varphi(\theta_\Delta(t)),
  \quad \forall k=0,1,2...
\]
Thus, it is convenient to assume that $\theta_\Delta$ mod $2\pi$
is a cyclic variable,
and the analysis is restricted to the range of $\theta_\Delta(0) \in [-\pi, \pi)$.

For the case of an odd PD characteristic\footnote{
There are examples of non odd PD characteristics, where \eqref{odd-change}
does not hold true (see, e.g.
BPSK Costas loop with sawtooth signals
\cite{LeonovKYY-2012-TCASII,LeonovKYY-2015-SIGPRO} and others).
}, system  \eqref{omega_delta_def} is not changed by the transformation
\begin{equation}\label{odd-change}
  \big(\omega_{\Delta}^{\text{free}},x(t),\theta_{\Delta}(t)) \rightarrow
  \big(-\omega_{\Delta}^{\text{free}},-x(t),-\theta_{\Delta}(t)).
\end{equation}
Property \eqref{odd-change} allows the analysis of system \eqref{final_system}
with only $\omega_\Delta^{\text{free}}>0$
and introduces the concept of \emph{frequency deviation}

\begin{center}{$|\omega_\Delta^{\text{free}}| = |\omega_1 - \omega_2^{\text{free}}|$.}\end{center}

\section{Locked state}
The locked states (also called steady states) of
the model in the signal's phase space must satisfy the following conditions:
\begin{itemize}
  \item the phase error $\theta_\Delta$ is constant, the frequency error $\dot \theta_\Delta$ is zero;
  \item the model in a locked state approaches the same locked state
  after small perturbations (of the VCO phase, input signal phase, and filter state).
\end{itemize}
The locally asymptotically stable equilibrium (stationary) points
of model \eqref{final_system}:
\begin{equation}
 \label{eq-points-def}
 \begin{aligned}
 & \theta_\Delta(t) \equiv \theta_{eq} + 2\pi k,\quad x(t) \equiv x_{eq},
 \end{aligned}
\end{equation}
are locked states, i.e. satisfy the above conditions\footnote{
It can be proved that if the filter is controllable and observable,
then only equilibria satisfy locked state conditions, i.e.
the filter state $x(t)$ must be constant in the locked state \cite{LeonovK-2014}.}.

\begin{figure}[ht]
 \centering
 \includegraphics[width=0.4\textwidth]{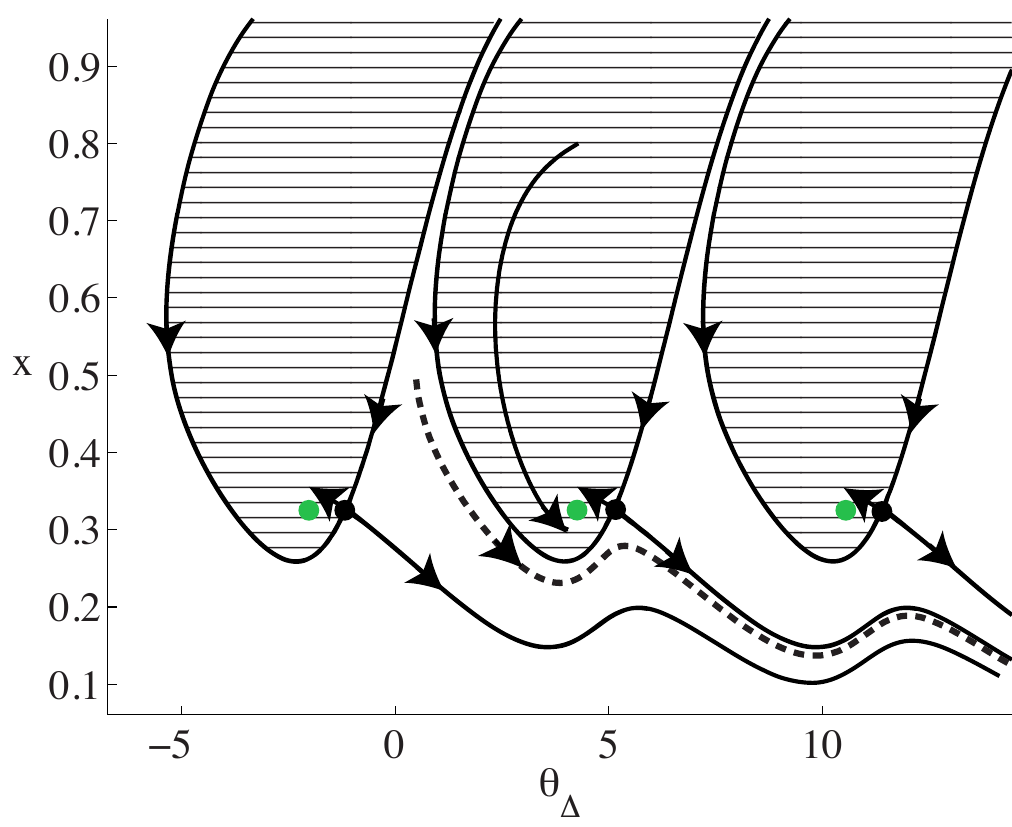}
 \caption{
  Phase portrait for $\omega_{\Delta}^{\text{free}}$ from the hold-in range.
  The system's evolving state over time traces a trajectory $(x(t),\theta_\Delta(t))$.
  Trajectories can not intersect.
  Unstable equilibrium points, such as saddles --- black dots,
  locally asymptotically stable equilibria --- green dots,
  any of which has its own basin of attraction (shaded domain)
  bounded by stable saddle separatrices (black trajectories going to the saddles).
  There are initial states and corresponding trajectories (see, e.g. dashed trajectory),
  which are not attracted to an equilibrium.
  }
 \label{hold-in-phase-portrait}
\end{figure}

Considering the case of a nonsingular matrix $A$
(i.e. the transfer function of the filter does not have zero poles),
the equilibria of \eqref{final_system} (stationary points)
are given by the equations
\begin{equation}\label{zeros1}
 \begin{aligned}
  & \varphi(\theta_{eq}) = \frac{\omega_{\Delta}^{\text{free}}}{L(-c^*A^{-1}b + h)}
   = \frac{\omega_{\Delta}^{\text{free}}}{LH(0)}, \\
  &  x_{eq} = -A^{-1}b \varphi(\theta_{eq}) =
  - A^{-1}b\frac{\omega_{\Delta}^{\text{free}}}{LH(0)}.
 \end{aligned}
\end{equation}
Thus, the equilibria can be considered as
a multiple-valued function of variable $\omega_\Delta^{\text{free}}$:
$\big(x_{eq}(\omega_{\Delta}^{\text{free}}),\theta_{eq}(\omega_{\Delta}^{\text{free}})\big)$.
From the boundedness of the PD characteristic $\varphi(\theta_{eq})$
it follows that there are no equilibria
for sufficiently large $|\omega_{\Delta}^{\text{free}}|$
(see Figs.~\ref{hold-in-phase-portrait} and \ref{hold-in-phase-portrait-noeq}).

\section{Engineering definitions of stability ranges}

The widely used engineering assumption (see Viterbi's pioneering writing
\cite[p.15]{Viterbi-1966}) is that
the zero input response of filter $\alpha_0(t,x_0)$  does not affect the synchronization of the loop.
This assumption allows the filter state $x(t)$ to be excluded
from the consideration and
a \emph{simplified mathematical model of PLL-based circuit in the signal's phase space}
to be obtained from \eqref{loop-filter-int} and \eqref{vco first}
(see, e.g. \cite[p.17, eq.2.20]{Viterbi-1966} for $h=0$
and \cite[p.41, eq.4-26]{Gardner-1966} for $\gamma \equiv 0$):
\begin{equation}\label{mathmodel-class-simple}
   \dot\theta_{\Delta}\!=\!
   \omega_{\Delta}^{\text{free}}-L\!\!\!\int_0^t\!\!\!\!\!\gamma(t - \tau)
   \varphi(\theta_{\Delta}(\tau)){\rm d}\tau
   -Lh \varphi(\theta_{\Delta}(t)).
\end{equation}
For an example of this one-dimensional integro-differential equation
the following intervals (\!\!\cite{Gardner-1966,Viterbi-1966}) are defined:
the hold-in range includes $|\omega_{\Delta}^{\text{free}}|$
such that model \eqref{mathmodel-class-simple}
has an equilibrium $\theta_{\Delta}(t) \equiv \theta_{eq}$,
which is locally stable
(local stability, i.e. for some initial phase error $\theta_{\Delta}(0)$);
the pull-in range includes $|\omega_{\Delta}^{\text{free}}|$
such that any solution of model \eqref{mathmodel-class-simple}
is attracted to one of the equilibria $\theta_{eq}$
(global stability, i.e. for any initial phase error $\theta_{\Delta}(0)$).
Thus, the block diagram of the loop in Fig.~\ref{costas-pll-sim-model}
is usually considered without initial data $x(0)$ and $\theta_{\Delta}(0)$
(see, e.g. \cite[p.17, fig.2.3]{Viterbi-1966}).

Viterbi \cite{Viterbi-1966} explains the above assumption for the stable matrix $A$,
but considers also various filters with marginally stable matrixes
(e.g. a filter -- perfect integrator, where $A=0$).
At the same time, even for a stable matrix $A$,
the initial filter state $x(0)$ and $\alpha_0(t,x_0)$
may affect the acquisition process and stability ranges
(see, e.g. corresponding examples for the classical PLL
\cite{KuznetsovKLNYY-2015-ISCAS} and Costas loops
\cite{KudryashovaKKLSYY-2014-ICINCO,KuznetsovKLNYY-2014-ICUMT-QPSK,KuznetsovKLSYY-2014-ICUMT-BPSK,BestKKLYY-2015-ACC}).

While the above assumption allows
introduction of the above one-dimensional stability sets,
defined only by $|\omega_{\Delta}^{\text{free}}|$,
for rigorous study the multi-dimensional stability domains have to considered,
taking into account $x(0)$,
and their relationships with the classical engineering ranges
have to be explained.
In \cite[p.187]{Gardner-2005-book} it is noted that the consideration of all state variables is of utmost importance
in the study of cycle slips and the \emph{lock-in} concept.
%Everything should be made as simple as possible, but not simpler

\smallskip

\section{Rigorous definitions of stability sets}
The rigorous mathematical definitions
of the hold-in, pull-in, and lock-in sets are now given
for the nonlinear mathematical model of PLL-based circuits in the signal's phase space
\eqref{final_system} and corresponding nontrivial examples are considered.

\subsection{Local stability and hold-in set}

We now consider the linearization\footnote{
Here it is assumed that
the PD characteristic $\varphi(\theta_\Delta)$ is smooth at the point $\theta_\Delta=\theta_{eq}$.
However, there are PLL-based circuits with nonsmooth or discontinuous PD characteristics
(see, e.g. the sawtooth PD characteristic for PLL \cite{Gardner-2005-book},
the model of QPSK Costas loop \cite{BestKLYY-2014-IJAC},
and some others \cite{biggio2014accurate,biggio2015efficient,bizzarri2012periodic}).
In such a case care has to be taken of the definition of solutions,
the linearization of the model
and the analysis of possible sliding solutions (see, e.g. \cite{GeligLY-1978}).}
of system \eqref{final_system} along
an equilibrium $(x_{eq},\theta_{eq})$.
Taking into account \eqref{zeros1} and $\varphi'(\theta) := d \varphi(\theta)/ d \theta$,
the linearized system is as follows:
\begin{equation}\label{linearized_system_2}
 \begin{aligned}
   \left(
     \begin{array}{c}
       \dot{x} \\
       \dot\theta_\Delta \\
     \end{array}
   \right)
   =
   \left(
   \begin{array}{cc}
     A & b\varphi'(\theta_{eq}) \\
     -Lc^* & -Lh\varphi'(\theta_{eq})
   \end{array}
   \right)
   \left(
     \begin{array}{c}
       x-x_{eq} \\
       \theta_\Delta - \theta_{eq} \\
     \end{array}
   \right)
 \end{aligned}
\end{equation}
The characteristic polynomial of linear system \eqref{final_system}
can be written (using the Schur complement, e.g. \cite{LeonovK-2014}) in the following form:
$\chi(s) = \big(-Lh\varphi'(\theta_{eq}) - s +Lc^*(A-sI)^{-1}b\varphi'(\theta_{eq})\big)\det(A-sI)$,
or can be expressed in terms of the filter's transfer function $H(s)=\frac{a(s)}{d(s)}$, where $a(s)$ and $d(s)$ are polynomials:
\begin{equation}\label{tfdenumerator}
  \begin{aligned}
    \chi(s) = -\big(sd(s)+a(s)L \varphi'(\theta_{eq})\big).
  \end{aligned}
\end{equation}
The characteristic polynomial corresponds to the denominator of the closed loop transfer function\footnote{
  Consideration of linearized model \eqref{linearized_system_2} allows to avoid
  the rigorous discussion of initial states $(x(0),\theta_\Delta(0))$
  related to the Laplace transformation and transfer functions \cite{LeonovK-2014}.
}.

To study the local stability of equilibria \eqref{zeros1},
it is necessary to check whether all the roots of the characteristic polynomial \eqref{tfdenumerator}
for the linearization of model \eqref{final_system} along the equilibria
(i.e. the poles of the closed loop transfer function)
have negative real parts.
For this purpose, at the stage of \emph{pre-design analysis}
when all parameters of the loop can be chosen precisely,
the Routh-Hurwitz criterion and its analogs
(see, e.g. Kharitonov's generalization \cite{Kharitonov-1978}
for interval polynomials) can be effectively applied.
At the stage of \emph{post-design analysis}
when only the input and VCO output are considered and the parameters
are known only approximately,
various frequency characteristics of the loop
(see, e.g. Nyquist and Bode plots)
and the continuation principle can be used (see, e.g, \cite{Gardner-2005-book,Best-2007}).

If the PD characteristic
is an odd function and hence $\varphi'(\theta_{eq})$ is an even function,
from \eqref{odd-change} we conclude that

1) there are symmetric equilibria:
$\big( x_{eq}(\omega_\Delta^{\text{free}}), \theta_{eq}(\omega_\Delta^{\text{free}}) \big)
=\big( -x_{eq}(-\omega_\Delta^{\text{free}}), -\theta_{eq}(-\omega_\Delta^{\text{free}}) \big)$,

2) these symmetric equilibria are simultaneously stable or unstable.

The same holds true for nonstationary trajectories.

\begin{definition}\label{def-hold}
  A set of all frequency deviations
  $|\omega_{\Delta}^{\text{free}}|$
  such that the mathematical model of the loop
  in the signal's phase space has a locally asymptotically stable equilibrium
  is called a hold-in set $\Omega_{\text{hold-in}}$.
\end{definition}

Thus, a value of frequency deviation belongs to the hold-in set
if the loop re-achieves locked state after small perturbations
of the filter's state, the phases and frequencies of VCO and the input signals.
This effect is also called \emph{steady-state stability}.
In addition, for a frequency deviation within the hold-in set,
the loop in a locked state tracks
small changes in input frequency,
i.e. achieves a new locked state (\emph{tracking process}).

In the literature the following explanations of the hold-in range
(sometimes also called a \emph{lock range} \cite[p.507]{PedersonM-2008-book}, \cite[p.10-2]{BakshiG-2009-book},
a \emph{synchronization range} \cite{Blanchard-1976}, a \emph{tracking range} \cite[p.49]{KiharaOE-2002})
can be found:
``{\it{The hold-in range is obtained by calculating the frequency where
the phase error is at its maximum}}''\cite[p.171]{brendel2013millimeter},
``{\it{The maximum frequency difference
before losing lock of the PLL system is called the hold-in range}}''\cite[p.258]{Kroupa-2012}.
The following example shows that these explanations may not be correct,
because for high-order filters
the hold-in ``range'' may have holes.

The following example shows that the hold-in set may not include
$\omega_\Delta^{\text{free}} = 0$.

\begin{example}[the hold-in set does not contain $\omega_\Delta^{\text{free}} = 0$]
Consider the classical PLL with the sinusoidal PD characteristic
$\varphi(\theta_\Delta) = \frac{1}{2}\sin(\theta_\Delta)$, VCO input gain $L=8$,
and the filter transfer function
 \begin{equation}
 H(s) = \frac{a(s)}{d(s)}= \frac{1+0.5s}{1+0.5s+0.5s^2}.
 \end{equation}
   From \eqref{zeros1} the following equation for equilibria is obtained:
  \begin{equation}
    \label{eq-points-1}
    \begin{aligned}
      & \frac{1}{2}\sin(\theta_{eq}) = \frac{1}{8}\omega_\Delta^{\text{free}}.
    \end{aligned}
  \end{equation}
\begin{figure}[ht]
 \centering
 \includegraphics[width=0.4\textwidth]{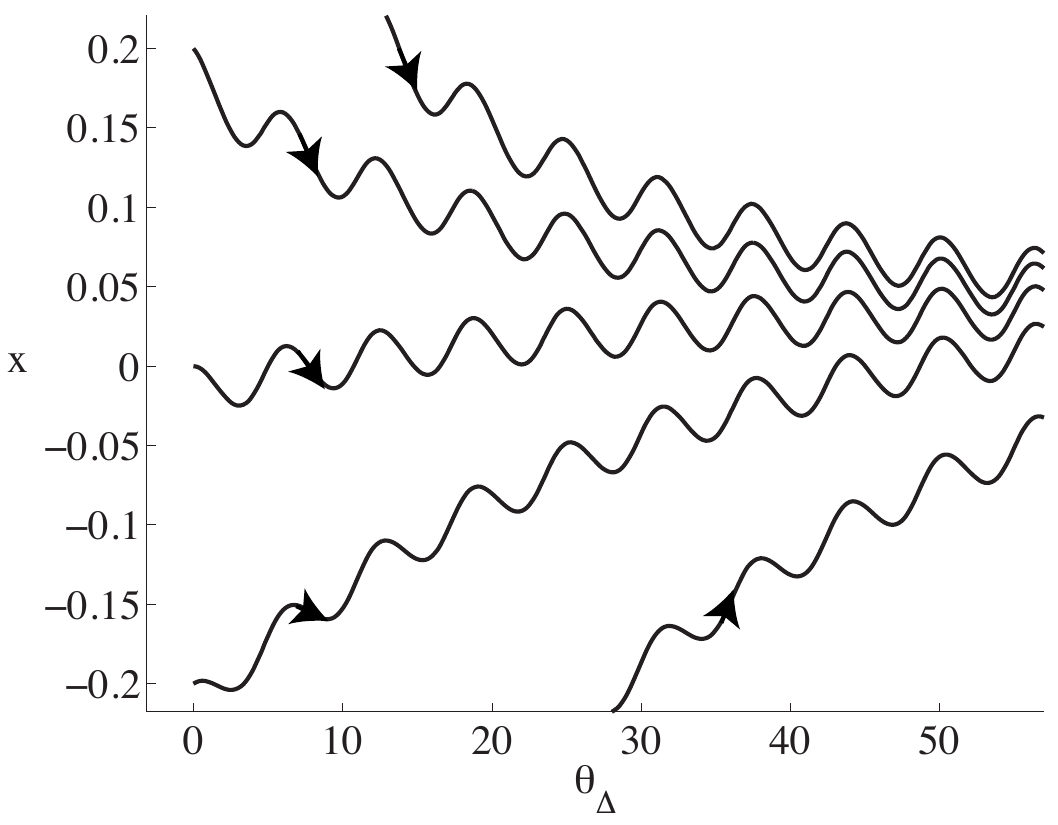}
 \caption{
 Phase portrait for $\omega_{\Delta}^{\text{free}}$ outside the hold-in range:
 there are no locally stable equilibria.}
 \label{hold-in-phase-portrait-noeq}
\end{figure}

  Applying the Routh-Hurwitz stability criterion\footnote{
  For a third-order polynomial $\chi(s) = a_3s^3 + a_2s^2 + a_1s + a_0$,
  all the roots have negative real parts and the corresponding linear system is asymptotically stable
  if $a_{1,2,3} > 0$ and  $a_2a_1 > a_3a_0$.
  For
  $\chi(s) = a_4s^4 + a_3s^3 + a_2s^2 + a_1s + a_0$,
  all the coefficients must satisfy  $a_{1,2,3,4} > 0$, and
  $a_3a_2 > a_4a_1$ and  $a_3a_2a_1 > a_4a_1^2 + a_3^2a_0$.
  }
  to the denominator of the closed loop transfer function \eqref{tfdenumerator}
  \begin{equation}
    \label{pol1}
    \begin{aligned}
      & s^3 + s^2 + s (2+4\cos(\theta_{eq})) + 8\cos(\theta_{eq}),
    \end{aligned}
  \end{equation}
  the following conditions are obtained:
  \begin{equation}
    \begin{aligned}
      & \cos(\theta_{eq}) >0, \ (2+4\cos(\theta_{eq})) > 0,\\
      & (2+4\cos(\theta_{eq})) > 8\cos(\theta_{eq}).
    \end{aligned}
  \end{equation}
Then $0< \cos(\theta_{eq}) < \frac{1}{2}$, and
for the locked state the steady-state phase error (i.e. corresponding to an equilibrium)
is obtained
\begin{equation}
    \label{theta-eqs-1}
    \begin{aligned}
      & \theta_{eq} \in (-\frac{\pi}{2},-\frac{\pi}{3})\cup
        (\frac{\pi}{3},\frac{\pi}{2}).
    \end{aligned}
  \end{equation}
Taking into account \eqref{eq-points-1}, \eqref{theta-eqs-1},
one obtains the hold-in set
 \begin{equation}
 \begin{aligned}
    & |\omega_\Delta^{\text{free}}| \in (2\sqrt{3},4).
 \end{aligned}
 \end{equation}
\end{example}
\smallskip

The next example shows that the hold-in set may not actually be a range
(i.e., an interval)
but a union of intervals, one of which may include
$\omega_\Delta^{\text{free}} = 0$.

\begin{example}[the hold-in set is a union of disjoint
intervals, one of which contains $\omega_\Delta^{\text{free}} = 0$]
\label{interval-union-example}
 Consider the classical PLL with the sinusoidal PD characteristic
 $\varphi(\theta_\Delta) = \frac{1}{2}\sin(\theta_\Delta)$, the VCO input gain $L=80$,
 and the filter transfer function
 \begin{equation}
   H(s) = \frac{1+0.25s+0.5s^2}{1+2s+2s^2+2s^3}.
 \end{equation}

From \eqref{zeros1} the following equation for the equilibria is obtained:
  \begin{equation}
  \label{eq-points-omega-2}
    \begin{aligned}
      & \frac{1}{2}\sin(\theta_{eq}) = \frac{1}{80}\omega_\Delta^{\text{free}}.
    \end{aligned}
  \end{equation}
An equilibrium is asymptotically stable if and only if all the roots
of polynomial \eqref{tfdenumerator}:
  \begin{equation}\label{pol2}
    \begin{aligned}
      & s(1+2s+2s^2+2s^3) + K(1+0.25s+0.5s^2) =\\
      & 2s^4 + 2s^3 + s^2(2+0.5K) + s(1+0.25K) + K,\\
      & K = L\varphi'(\theta_{eq})= 40\cos(\theta_{eq})
    \end{aligned}
  \end{equation}
have negative real parts.
Using the Routh-Hurwitz criterion, we obtain
\begin{equation}
  \begin{aligned}
    & 2+0.5K >0,\ 1+0.25K > 0,\ K > 0,\\
    & 2(2+0.5K) > 2(1+0.25K),\\
    & 2(2+0.5K)(1+0.25K) > 2(1+0.25K)^2+2^2K.
  \end{aligned}
\end{equation}
From these inequalities we have
\begin{equation}
\label{theta-eqs-2}
  \begin{aligned}
      & K = 40\cos(\theta_{eq}) \in (0,12 - 8\sqrt{2})\cup(12+8\sqrt{2},\infty),\\
      & \theta_{eq} \in (-\frac{\pi}{2},-1.5536) \cup (-0.9486,0.9486) \cup (1.5536,\frac{\pi}{2}).
    \end{aligned}
  \end{equation}
Note that for other values of $\theta_{eq}$
at least one root of the polynomial \eqref{pol2}
has a positive real part, making the corresponding equilibrium unstable.
Combining \eqref{eq-points-omega-2} and \eqref{theta-eqs-2},
we obtain the hold-in set
\begin{equation}
 \begin{aligned}
 & |\omega_\Delta^{\text{free}}| \in [0,32.5) \cup (39.9942,40).
 \end{aligned}
\end{equation}
Note that in this case, for the values of the VCO input gain $L > 24+16\sqrt{2}$
the hold-in set is always a union of disjoint intervals.
For $L=80$ the simulation results of transition process
in Simulink model\,\footnote{Following
the above classical consideration,
the filter is often represented in MatLab Simulink as the block \emph{Transfer Fcn}
with zero initial state
(see, e.g. \cite{BrigatiFMPP-2001,NicolleTMOJ-2007,Zucchelli-2007,KoivoE-2009,KaaldLHS-2009}).
It is also related to the fact that the transfer function (from $\varphi$ to $g$)
of linear system \eqref{loop-filter}
is defined by the Laplace transformation for zero initial data $x(0) \equiv 0$.
In Fig.~\ref{simulink-model} we use the block \emph{Transfer Fcn (with initial states)}
to take into account the initial filter state $x(0)$;
the initial phase error $\theta_\Delta(0)$
can be taken into account by the property \emph{initial data} of the \emph{Intergator} blocks.
Note that the corresponding initial states in SPICE
(e.g. capacitor's initial charge) % and inductor's initial currents
are zero by default
but can be changed manually \cite{BianchiKLYY-2015}.}
in Fig.~\ref{simulink-model} are shown in Figs.~\ref{sim-1}--\ref{sim-3}
for the initial data $(x(0)=(0;0;0.9990), \theta_\Delta(0) = 1.5585)$
and various $\omega_\Delta^{\text{free}}$.
\begin{figure}[ht]
\centering
  \includegraphics[width=0.48\textwidth]{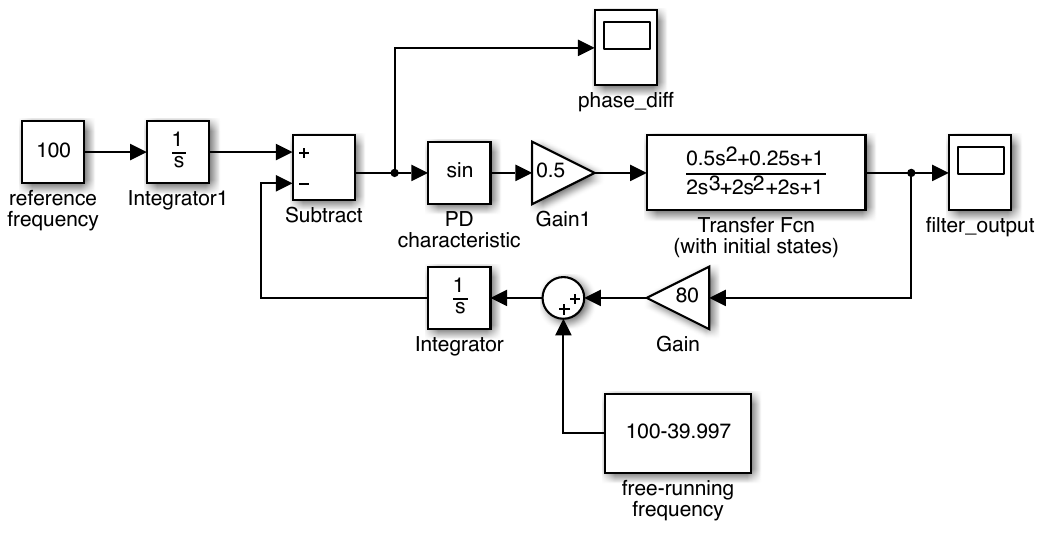}
  \caption{MatLab Simulink: the signal's phase space model of the classical PLL}
  \label{simulink-model}
\end{figure}
\end{example}

Related discussion on the frequency responses of loop
with high-order filters can be found in \cite[p.34-38, 52-56]{Gardner-2005-book}.

\begin{remark}
  For the first order filters, the set $\Omega_{\text{hold-in}}$
  is an interval $|\omega_\Delta^{\text{free}}| < \omega_h$.
  For higher order filters, the set $\Omega_{\text{hold-in}}$
  may be more complex.
  Thus, from an engineering point of view,
  it is reasonable to require that
  $\omega_\Delta^{\text{free}} = 0$ belongs to the hold-in set
  and to define a hold-in range as
  the largest interval $[0,\omega_{h})$ from the hold-in set
  \[
    [0,\omega_{h}) \subset \Omega_{\text{hold-in}}
  \]
  such that a certain stable equilibrium varies continuously
  when $\omega_\Delta^{\text{free}}$ is changed within the range\footnote{\label{fnote}
  In general (when the stable equilibria coexist and some of them may appear or disappear),
  the stable equilibria can be considered as
  a multiple-valued function of variable $\omega_\Delta^{\text{free}}$,
  in which case the existence of its continuous singlevalue branch
  for $|\omega_\Delta^{\text{free}}| \in [0,\omega_{h})$ is required.}.
  Here $\omega_{h}$ is called
  a \emph{hold-in frequency} (see \cite[p.38]{Gardner-1966}).
\end{remark}

\begin{remark}
In the general case when there is no symmetry with respect to
$\omega_{\Delta}^{\text{free}}$ (see \eqref{odd-change})
the hold-in set need not be symmetric
and the set $\omega_{\Delta}^{\text{free}} \in \Omega_{\text{hold-in}}$
must be considered  in Definition~\ref{def-hold}.
\end{remark}

\begin{figure}[ht]
 \centering
 \includegraphics[width=0.24\textwidth]{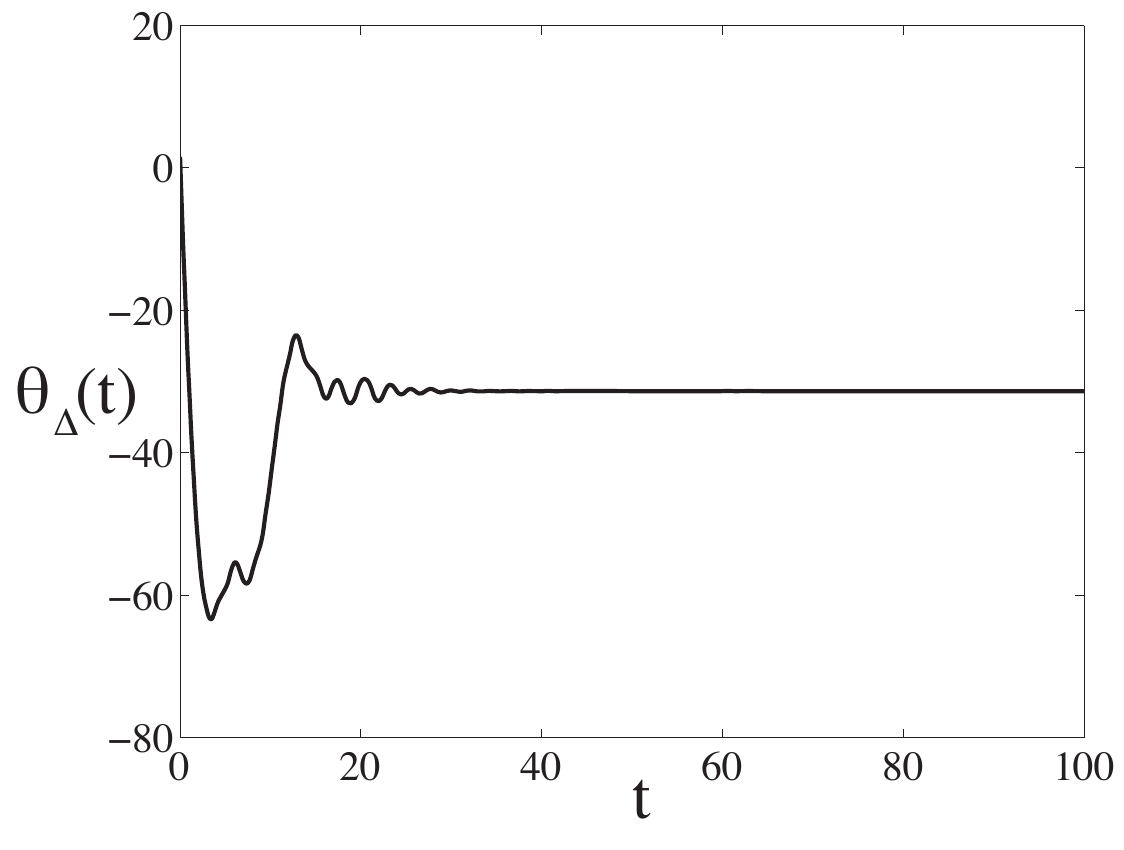}
 \includegraphics[width=0.24\textwidth]{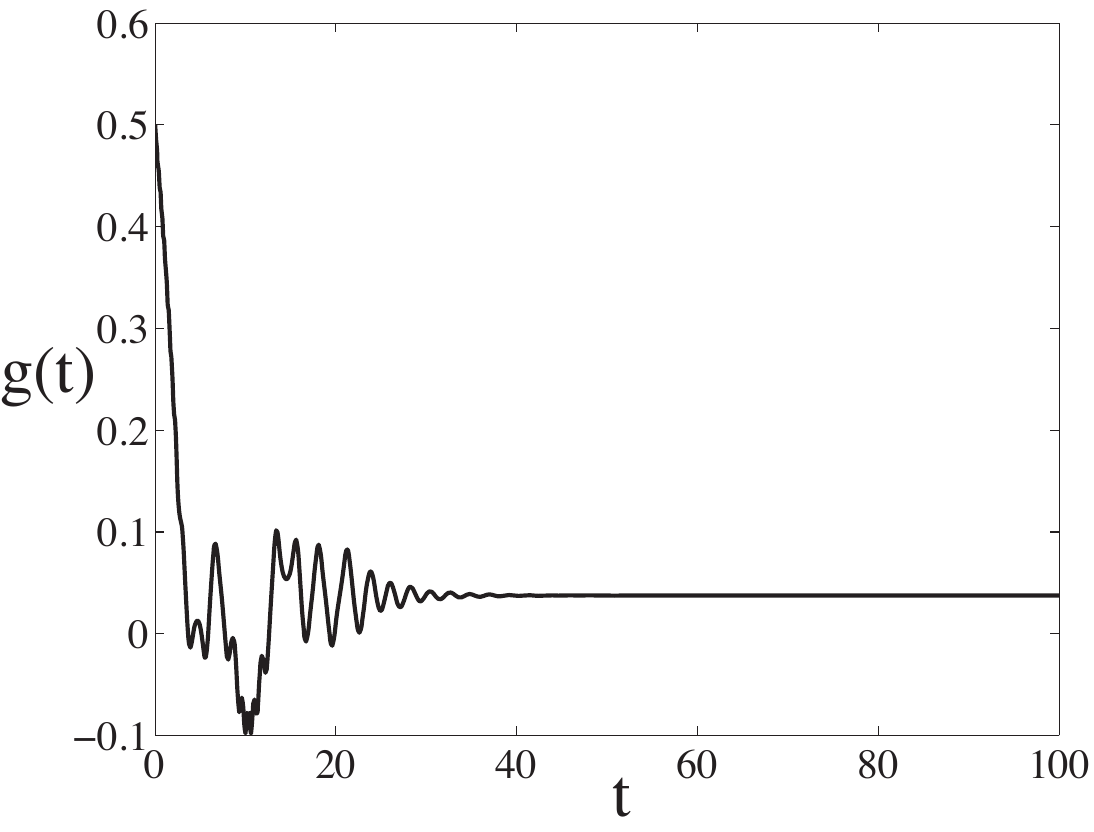}
 \caption{$\omega_\Delta^{\text{free}} = 3$: stable locked state exists.}
 \label{sim-1}
\end{figure}
\begin{figure}[ht]
 \centering
 \includegraphics[width=0.24\textwidth]{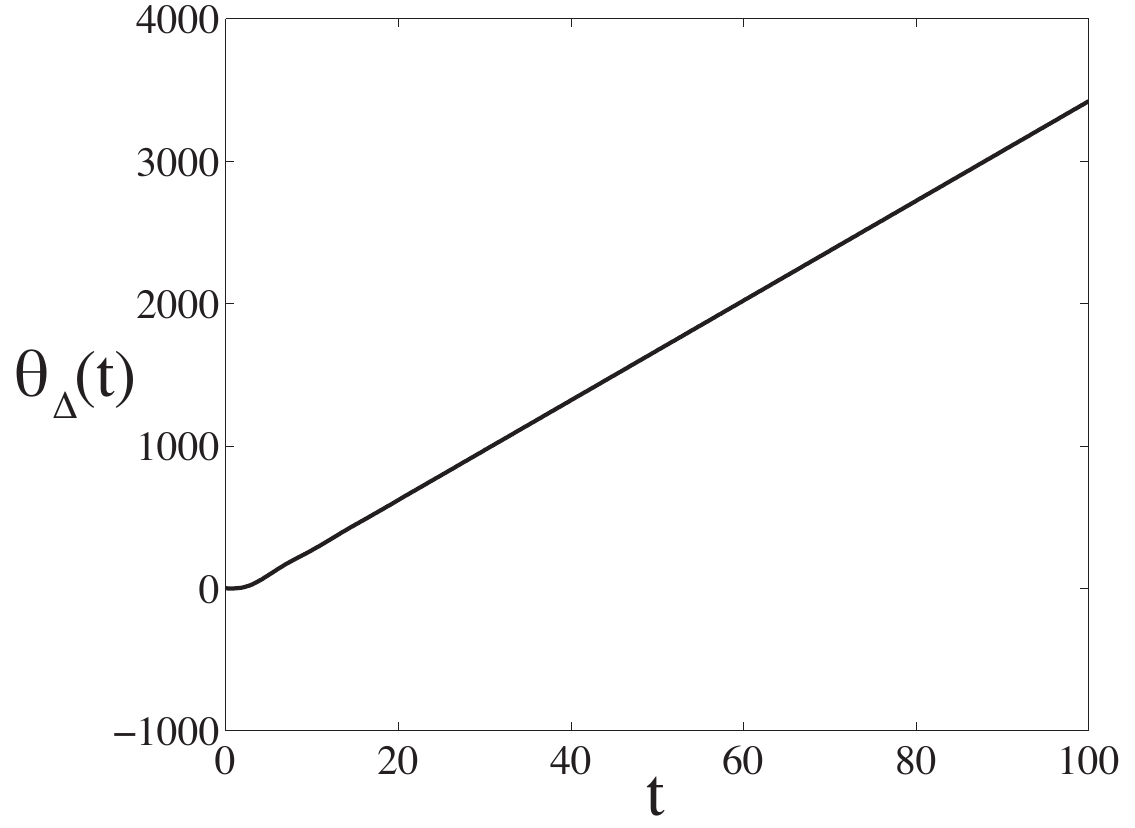}
 \includegraphics[width=0.24\textwidth]{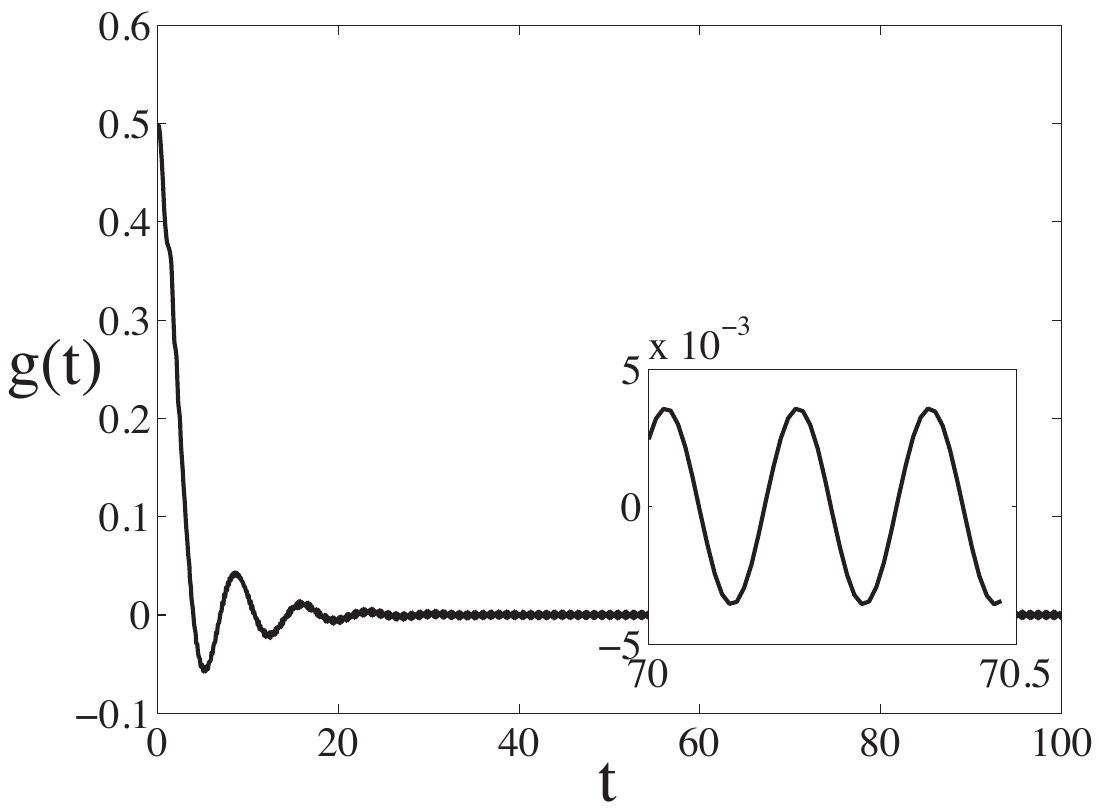}
 \caption{$\omega_\Delta^{\text{free}} = 35$: there are no locked states (see also Fig.~\ref{hold-in-phase-portrait-noeq}).}
 \label{sim-2}
\end{figure}
\begin{figure}[ht]
 \centering
 \includegraphics[width=0.24\textwidth]{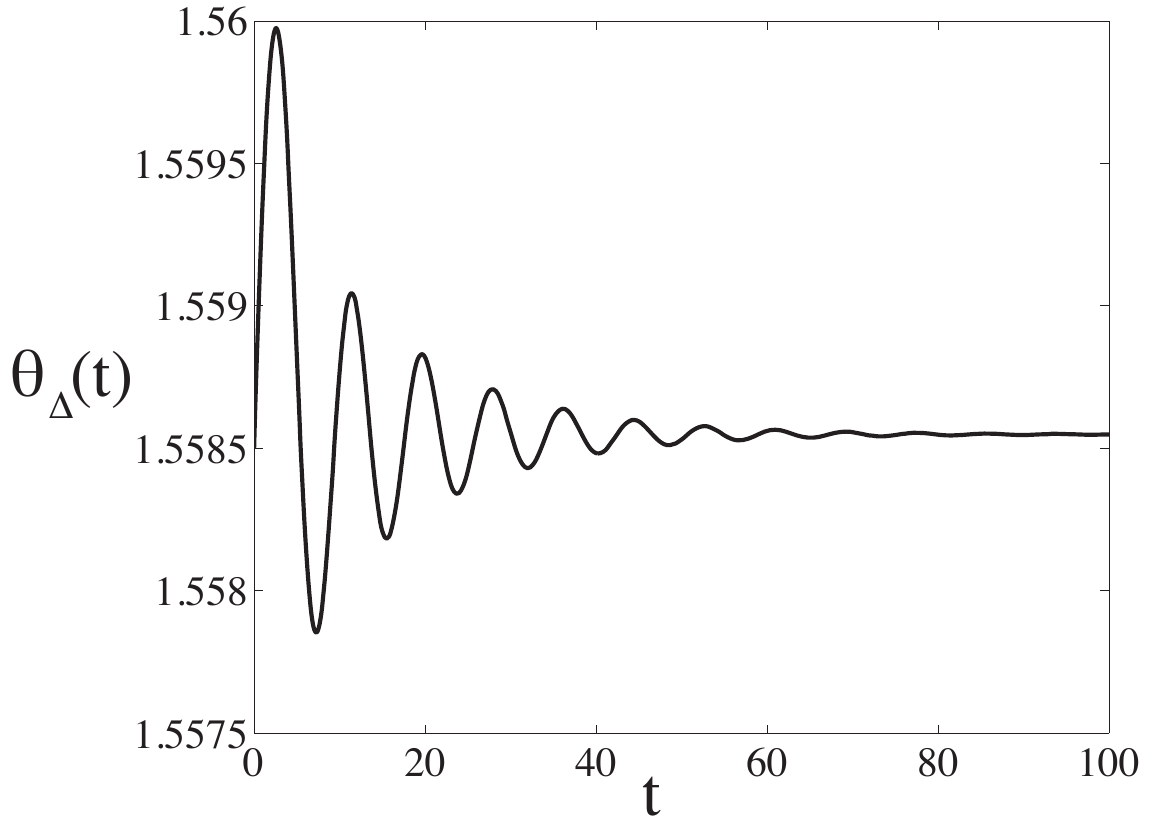}
 \includegraphics[width=0.24\textwidth]{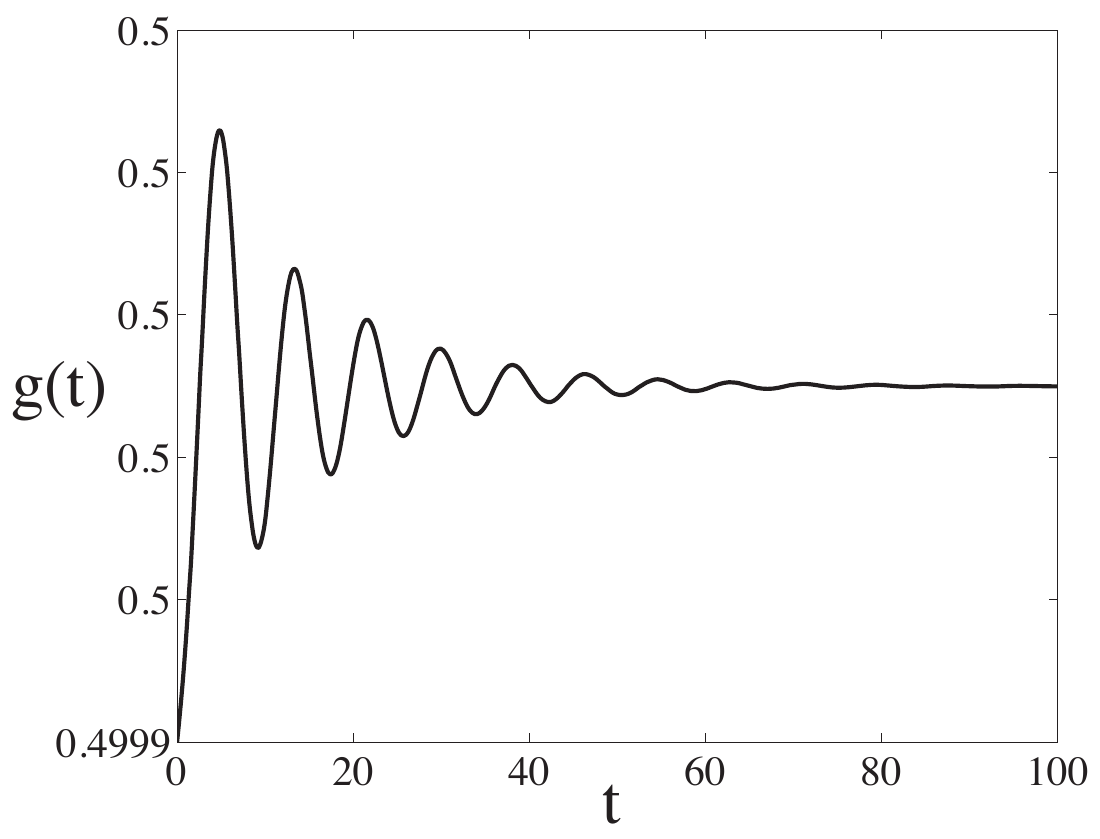}
 \caption{$\omega_\Delta^{\text{free}} = 39.997$: stable locked state exists.}
 \label{sim-3}
\end{figure}

\subsection{Global stability (stability in the large) and pull-in set}

Assume that the loop power supply is initially switched off and then at $t = 0$ the power is switched on,
and assume  that the initial frequency difference is sufficiently large.
The loop may not lock within one beat note,
but the VCO frequency will be slowly tuned
toward the reference frequency (\emph{acquisition process}).
This effect is also called a \emph{transient stability}.
The pull-in range is used to name such frequency deviations
that make the acquisition process possible
(see, e.g. explanations in
\cite[p.40]{Gardner-1966}, \cite[p.61]{Best-2007}).

To define a \emph{pull-in range}
(called also a \emph{capture range} \cite{Talbot-2012-book},
an \emph{acquisition range} \cite[p.253]{Blanchard-1976})
rigorously, consider first an important definition from stability theory.

\begin{figure}[ht]
 \centering
 \includegraphics[width=0.36\textwidth]{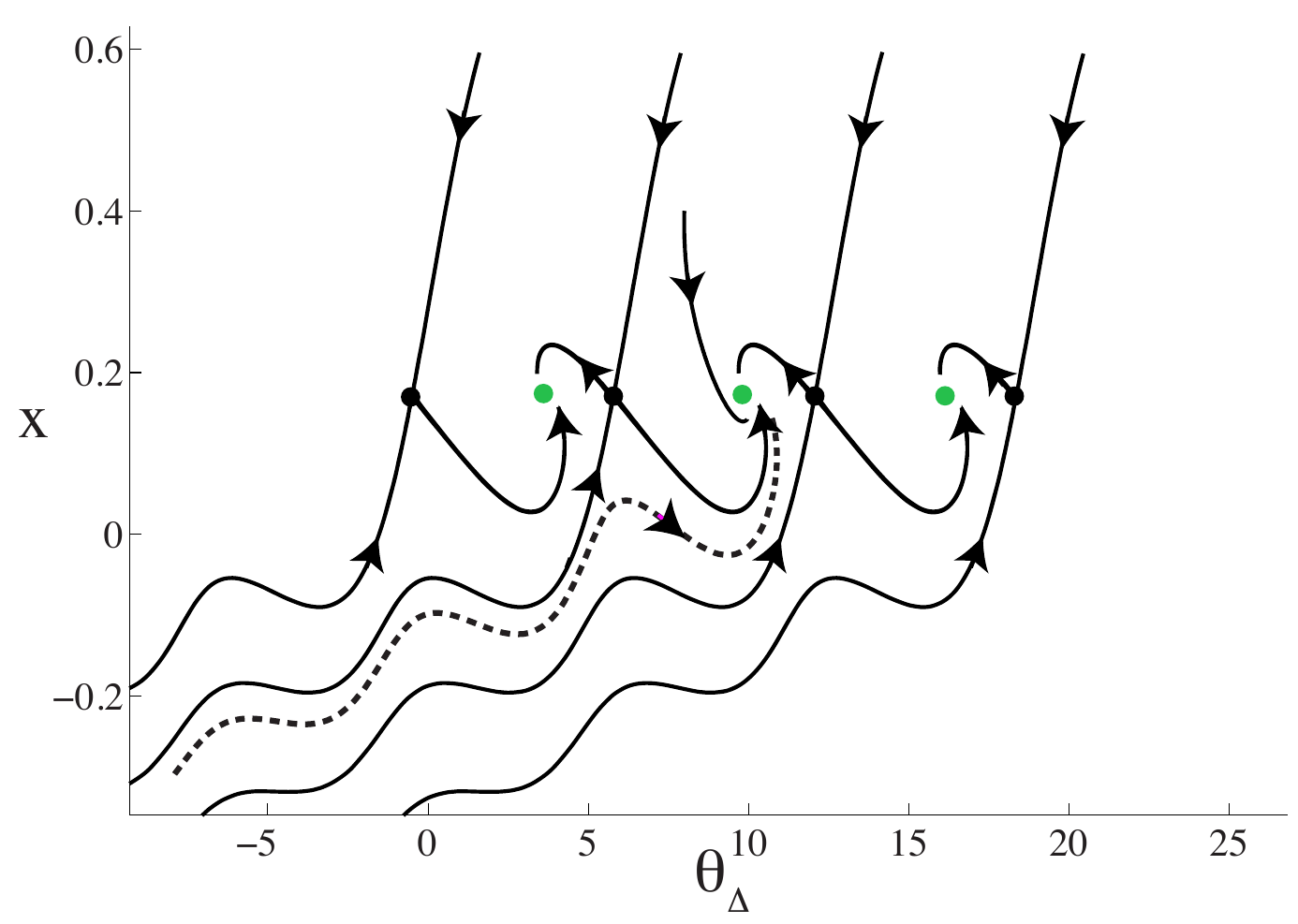}
 \caption{ Phase portrait for $\omega_{\Delta}^{\text{free}}$
  from the pull-in range:
  any trajectory is attracted to an equilibrium
  (equilibria: green--stable and black--unstable circles);
  for a sufficiently large initial state of the filter,
  cycle slipping is possible (see, e.g. dashed trajectory).}
 \label{fig-pullin}
\end{figure}

\begin{definition}
  If for a certain $\omega_{\Delta}^{\text{free}}$
  any trajectory of system \eqref{final_system} tends to an equilibrium,
  then the system with such $\omega_{\Delta}^{\text{free}}$
  is called globally asymptotically stable
  (see Fig.~\ref{fig-pullin}).
 \end{definition}

We now consider a possible rigorous definition.

\begin{definition} \label{def-pull}
 A set of all frequency deviations $|\omega_\Delta^{\text{free}}|$
 such that the mathematical model of the loop in the signal's phase space
 is globally asymptotically stable
 is called a pull-in set $\Omega_{\text{pull-in}}$.
\end{definition}

\begin{remark}
  In the general case when there is no symmetry with respect to
  $\omega_{\Delta}^{\text{free}}$
  the set $\omega_{\Delta}^{\text{free}} \in \Omega_{\text{pull-in}}$
  has to be considered in Definition~\ref{def-pull}.
\end{remark}

\begin{remark}
  The pull-in set is a subset of the hold-in set:
  $\Omega_{\text{pull-in}} \subset \Omega_{\text{hold-in}}$,
  and need not be an interval.
  From an engineering point of view,
  it is reasonable to require that
  $\omega_\Delta^{\text{free}} = 0$ belongs to the pull-in set
  and to define a pull-in range as
  the largest interval $[0,\omega_{p})$ from the pull-in set:
  \[
    [0,\omega_{p}) \subset \Omega_{\text{pull-in}},
  \]
  where $\omega_{p}$ is called a \emph{pull-in frequency} (see \cite[p.40]{Gardner-1966}).
\end{remark}

\begin{remark}
   If all possible states of the filter are bounded:
   \[
     x \in X_{\text{real}}\ (\text{e.g.\ } X_{\text{real}} =
     \{x: c_{\text{min}} <|x| < c_{\text{max}}\}),
   \]
   by the design of the circuit
   (e.g. capacitors have limited maximum and minimum charges,
   the VCO frequency is limited etc.),
   then in the definition of pull-in set
   it is reasonable to require that only
   solutions with $x(0) \in X_{\text{real}}$
   tend to the stationary set.
   Trajectories, with initial data outside of
   the domain defined by $x(0) \in X_{\text{real}}$
   (here the initial phase error $\theta_\Delta(0)$ can take any value),
   need not tend to the stationary set.
\end{remark}

For the model without filter (i.e. $H(s)=const$)
the pull-in set coincides with the hold-in set.
The pull-in set of PLL-based circuits with first-order filters
can be estimated using phase plane analysis methods \cite{Tricomi-1933,AndronovVKh-1937}, %,BarbashinT-1969
but in general its rigorous study is a challenging task
\cite{Viterbi-1966,Stensby-1997,Margaris-2004,KudrewiczW-2007,PinheiroP-2014}.

For the case of the passive lead-lag filter
$H(s) = \frac{1+s \tau_2}{1+s(\tau_1 + \tau_2)}$,
a recent work \cite[p.123]{Margaris-2004} notes that
``{\it the determination of the width of the capture range together
with the interpretation of the capture effect in the second order type-I loops
have always been an attractive theoretical problem.
This problem has not yet been provided with a satisfactory solution}''.
At the same time in \cite{Shakhtarin-1969,Belyustina-1970-eng,LeonovK-2013-IJBC,LeonovK-2014} %Gubar-1961
it is shown that the basin of attraction of the stationary set may be bounded
(e.g. by a semistable periodic trajectory, which
may appear as the result of collision of unstable and stable periodic solutions),
and corresponding analytical estimations and bifurcation diagram are given.

\begin{figure}[h]
  \centering
  \includegraphics[width=0.72\linewidth]{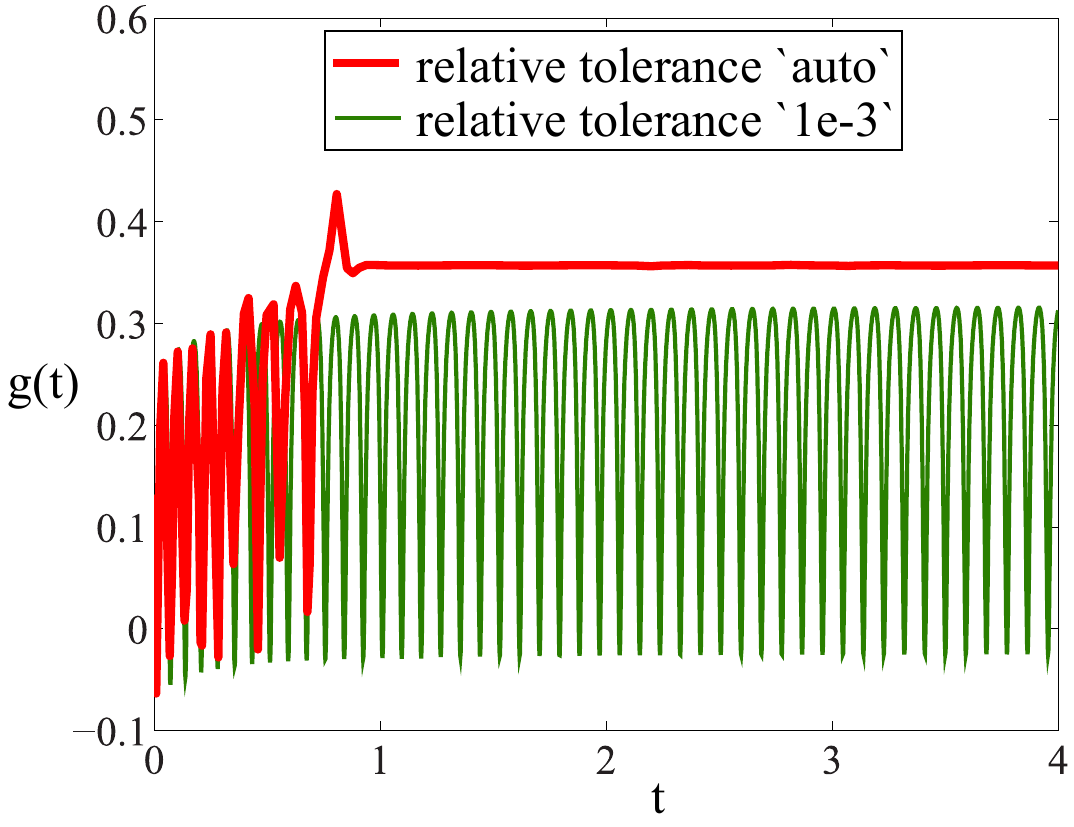}
  \caption{
  Simulation of two-phase PLL described by Fig.~\ref{simulink-model} or model \eqref{final_system} \cite{BianchiKLYY-2015}:
  $\tau_1\!=\!0.0448$, $\tau_2\!=\!0.0185$,
  $A\!=\!-\frac{1}{\tau_1+\tau_2}$, $b\!=\!1 - \frac{\tau_2}{\tau_1+\tau_2}$,
  $c\!=\!\frac{1}{\tau_1+\tau_2}$, $h\!=\!\frac{\tau_2}{\tau_1+\tau_2}$; $\varphi(\theta_\Delta)=\frac{1}{2}\sin(\theta_\Delta)$;
  $\omega_1\!=\!10000$, $\omega_2^{\text{free}}\!=\!10000 - 178.9$, $L\!=\!500$.
  Filter output $g(t)$ for the initial data $x_0\!=\!0.1318, \theta_{\Delta}(0)\!=\!0$
  obtained for default ``auto'' relative tolerance (red) --- acquires lock,
  relative tolerance set to ``1e-3''(green) --- does not acquire lock.}
  \label{pll_hidden}
\end{figure}

Note that in this case a numerical simulation may give wrong estimates
and should be used very carefully.
For example, in \cite{BianchiKLYY-2015}
the SIMetrics SPICE model for a two-phase PLL with a lead-lag filter
gives two essentially different results of simulation
with default ``auto'' sampling step (acquires lock) and minimum sampling step set to $1m$
(does not acquire lock --- such behaviour agrees with the theoretical analysis).
The same problems
are also observed in MatLab Simulink \cite{KuznetsovLYY-2014-IFAC,KuznetsovKLNYY-2015-ISCAS,BestKKLYY-2015-ACC},
see, e.g. Fig.~\ref{pll_hidden}.
These examples demonstrate the difficulties of numerical search of so-called \emph{hidden oscillations}
\cite{LeonovK-2013-IJBC,KuznetsovL-2014-IFACWC,LeonovKM-2015-EPJST},
whose basin of attraction does not overlap with the neighborhood of an equilibrium point,
and thus may be difficult to find numerically\footnote{In \cite{LauvdalMF-1997}
the crash of aircraft YF-22 Boeing in April 1992,
caused by the difficulties of rigorous analysis and design of nonlinear control systems with saturation,
is discussed and the conclusion is made
that since stability in simulations does not imply
stability of the physical control system,
stronger theoretical understanding is required
(see, e.g. similar problem  with the simulation of PLL in Fig.~\ref{pll_hidden}).
These difficulties in part are related to well-known Aizerman's and Kalman's conjectures on
the global stability of nonlinear control systems,
which are valid from the standpoint of
simplified analysis by the linearization, harmonic balance, and describing function methods
(note that all these methods are also widely used to the analysis of nonlinear oscillators used in VCO \cite{Margaris-2004,Suarez-2009}).
However the counterexamples (multistable high-order nonlinear systems where
the only equilibrium, which is stable, coexists with a hidden periodic oscillation) can be constructed
to these conjectures \cite{LeonovK-2013-IJBC,HeathCS-2015}.}.
In this case the observation of one or another stable solution may depend on the initial data
and integration step.

\noindent S.~Goldman, who has worked at Texas Instruments over 20 years, notes that
PLLs are used as pipe cleaners for breaking simulation tools \cite[p.XIII]{Goldman-2007-book}.

While PLL-based circuits are nonlinear control systems and
for their nonlocal analysis it is essential to apply
the classical stability criteria, which are developed in control theory,
however their direct application to analysis of the PLL-based models is often impossible,
because such criteria are usually not adapted for
the cylindrical phase space\footnote{For example,
in the classical Krasovskii--LaSalle principle on global stability
the Lyapunov function has to be radially unbounded
(e.g. $V(x,\theta_\Delta) \to +\infty$ as $||(x,\theta_\Delta)|| \to +\infty$).
While for the application of this principle to the analysis of phase synchronization systems
there are usually used Lyapunov functions periodic in $\theta_\Delta$
(e.g. $V(x,\theta_\Delta)$ in Remark~\ref{remarkPI} is bounded
for any $||(0,\theta_\Delta)|| \to +\infty$),
and the discussion of this gap is often omitted
(see, e.g. patent \cite{Abramovitch-2004} and works \cite{Bakaev-1963,Abramovitch-1990,Abramovitch-2003}).
Rigorous discussion can be found, e.g. in \cite{GeligLY-1978,LeonovK-2014}.
};
in the tutorial {\it Phase Locked Loops: a  Control Centric Tutorial}
\cite{Abramovitch-2002},
presented at {\it the American Control Conference 2002}, it was said that
``{\it The general theory of PLLs and ideas on how to make
them even more useful seems
to cross into the controls literature only rarely}''.

At the same time the corresponding modifications of classical stability criteria
for the nonlinear analysis of control systems in cylindrical phase space
were well developed in the second half of the 20th century,
see, e.g. \cite{GeligLY-1978,LeonovRS-1992,LeonovPS-1996,LeonovBSh-1996}.
A comprehensive discussion and the current state of the art can be found in \cite{LeonovK-2014}.
One reason why these works have remained
almost unnoticed by the contemporary engineering community
may be that they were written in
the language of control theory and the theory of dynamical systems,
and, thus, may not be well adapted to the terms and objects used
in the engineering practice of phase-locked loops.
Another possible reason, as noted in \cite[p.1]{Tranter-2010-book},
is that the nonlinear analysis techniques are well beyond the scope
of most undergraduate courses in communication theory and circuits design.
Note that for the application of various stability criteria
it is often necessary to represent system \eqref{final_system} in the Lur'e form:
\begin{equation}\label{Lurie}
 \begin{aligned} &
 \begin{aligned}
   \left(\!\!\!
     \begin{array}{c}
       \dot{\bar{x}} \\
       \dot\theta_\Delta\\
     \end{array}
   \!\!\!\right)
   =\left(\!\!\!
      \begin{array}{cc}
        A & 0 \\
        -Lc^* & 0 \\
      \end{array}
    \!\!\!\right)
   \left(\!\!\!
     \begin{array}{c}
       \bar{x} \\
       \theta_\Delta\\
     \end{array}
   \!\!\!\right)
   +
   \left(\!\!\!
     \begin{array}{c}
       b \\
       -Lh \\
     \end{array}
   \!\!\!\right)
   \bar{\varphi}(\theta_\Delta)
 \end{aligned}, \\
 \end{aligned}
\end{equation}
where
\[
 \begin{aligned}
  & \bar{x} = x-x_{eq} = x+A^{-1}b\varphi(\theta_{eq}),\
  \bar{\varphi}(\theta_\Delta) = \varphi(\theta_\Delta)-\varphi(\theta_{eq}), \\
  & \varphi(\theta_{eq}) = \omega_{\Delta}^{\text{free}}L^{-1}(-c^*A^{-1}b + h)^{-1}.
 \end{aligned}
\]

See also discussion of some nonlinear methods for the analysis of PLL-based models
in recent books \cite{SuarezQ-2003,Margaris-2004,KudrewiczW-2007,Suarez-2009}.

\subsection{Cycle slips and lock-in range}

Let us rigorously define \emph{cycle slipping} in the phase space of system \eqref{final_system}.
\begin{definition}\label{def-cs}
If
 \begin{equation}\label{eq-cs}
 \begin{aligned}
 & \limsup\limits_{t\to+\infty} |\theta_\Delta(0) - \theta_\Delta(t)| > 2\pi,
 \end{aligned}
 \end{equation}
it is then said that cycle slipping occurs
(see, e.g. dashed trajectory in Fig.~\ref{fig-pullin}).
\end{definition}

Here, sometimes, instead of the limit of the difference,
 the maximum of the difference is considered
(see, e.g. \cite[p.131]{Stensby-1997}).

\noindent{\bf Definition \ref{def-cs}'}
{\it
If
\begin{equation}\label{eq-cs-sup}
\begin{aligned}
   & \sup\limits_{t>0} |\theta_\Delta(0) - \theta_\Delta(t)| > 2\pi,
\end{aligned}
\end{equation}
it is then said that cycle slipping has occurred.
}

Note that, in general, Definition~\ref{def-cs}' need not mean that finally
(after acquisition) condition \eqref{eq-cs} can not be fulfilled.

Sometimes, the number of cycle slips is of interest.

\begin{definition}
If
\begin{equation}
\begin{aligned}
 & 2k\pi < \limsup\limits_{t\to\infty}
           |\theta_\Delta(0) - \theta_\Delta(t)| < 2(k+1)\pi,
\end{aligned}
\end{equation}
it is then said that $k$ cycle slips occurred.
\end{definition}

A numerical study of cycle slipping in classical PLL can be found in
\cite{AscheidM-1982}.
Analytical tools for estimating the number of cycle slips
depending on the parameters of the loop
can be found, e.g. in \cite{ErshovaL-1983,LeonovRS-1992,LeonovK-2014}.

The concepts of \emph{lock-in frequency} and \emph{lock-in range}
(called also a \emph{lock range}\cite[p.256]{Yeo-2010-book}, a  \emph{seize range} \cite[p.138]{Egan-2007-book}),
were intended to describe the set of frequency deviations for which
the loop can acquire lock within one beat without cycle slipping.
In \cite[p.40]{Gardner-1966} the following definition was introduced:
``{\it{If, for some reason,  the frequency difference between input and VCO
is less than the loop bandwidth, the loop will lock up almost instantaneously
without slipping cycles. The maximum frequency difference for which
this fast acquisition is possible is called the lock-in frequency}}''.

\begin{figure*}[ht]
 \centering
 \includegraphics[width=\textwidth]{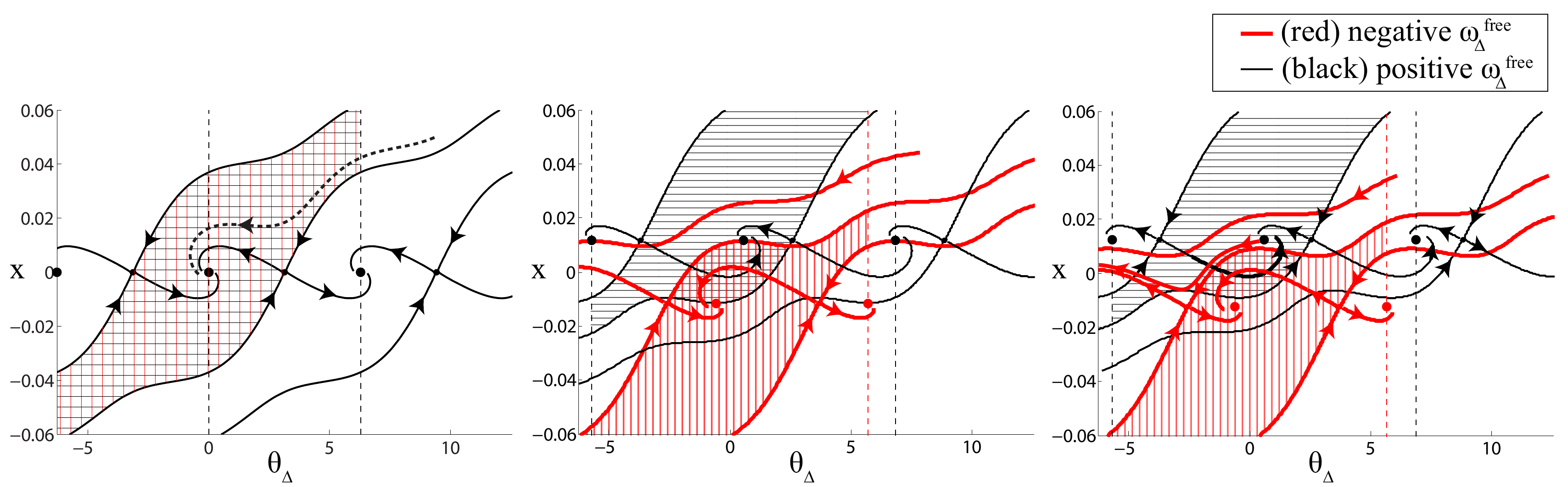}
 \caption{
 Phase portraits for the classical PLL with the following parameters:
 $H(s)= \frac{1+s\tau_2}{1+s(\tau_1 + \tau_2)}$,
 $\tau_1 = 4.48\cdot10^{-2}$,
 $\tau_2 = 1.85\cdot10^{-2}$,
 $L=250$, $\varphi(\theta_\Delta)= \frac{1}{2}\sin(\theta_\Delta)$, and various frequency deviations.
 Black color is for the system with positive $\omega_\Delta^{\text{free}}=|\widetilde{\omega}|$.
 Red is for the system with negative $\omega_\Delta^{\text{free}}=-|\widetilde{\omega}|$.
 Equilibria (dots), separatrices pass in and out of the saddles,
 local lock-in domains are shaded
 (upper black horizontal lines is for $\omega_\Delta^{\text{free}}>0$,
 lower red vertical lines is for $\omega_\Delta^{\text{free}}<0$).
 Left subfig: $\omega_\Delta^{\text{free}} = 0$;
 middle subfig: $\omega_\Delta^{\text{free}} = \pm 65$;
 right subfig: $\omega_\Delta^{\text{free}} = \pm 68$.
 }. \label{lockin3def}
\end{figure*}

However, in general, even for zero frequency deviation ($\omega_\Delta^{\text{free}}=0$)
and a sufficiently large initial state of filter ($x(0)$),
cycle slipping may take place (see, e.g. dashed trajectory in Fig.~\ref{lockin3def},\,left).
Thus, considering of all state variables is of utmost importance
for the cycle slip analysis and, therefore,
the concept \emph{lock-in frequency}
lacks rigor for classical simplified model \eqref{mathmodel-class-simple}
because it does not take into account the initial state of the filter.
The above definition of the lock-in frequency and corresponding definition of the lock-in range
were subsequently in various engineering publications
(see, e.g.
\cite[p.34-35]{Best-1984},\cite[p.161]{Wolaver-1991},\cite[p.612]{HsiehH-1996},\cite[p.532]{Irwin-1997},\cite[p.25]{CraninckxS-1998-book},
\cite[p.49]{KiharaOE-2002},\cite[p.4]{Abramovitch-2002},\cite[p.24]{DeMuerS-2003-book},\cite[p.749]{Dyer-2004-book},\cite[p.56]{Shu-2005},
\cite[p.112]{Goldman-2007-book},\cite[p.61]{Best-2007},\cite[p.138]{Egan-2007-book},\cite[p.576]{Baker-2011},\cite[p.258]{Kroupa-2012}).

The loop model \eqref{final_system}
has a subdomain of the phase space, where trajectories do not slip cycles
(called a lock-in domain), for each value of $\omega_\Delta^{\text{free}}$.
The lock-in domain is the union of local lock-in domains,
each of which corresponds to one of the equilibria and has
its own shape
(see, e.g. shaded domain in Fig.~\ref{lockin3def},\,left defined by corresponding separatrices).
The shape of the lock-in domain significantly varies depending
on $\omega_{\Delta}^{\text{free}}$.
In \cite[p.50]{Viterbi-1966}) a lock-in domain is called a \emph{frequency lock}.
Some writers (e.g. \cite[p.132]{Stensby-1997},\cite[p.355]{Meyer-2004-book})
use the concept \emph{lock-in range} to denote a \emph{lock-in domain}.

In general, taking into account nonuniform behavior of the lock-in domain shape,
Gardner wrote {\it ``There is no natural way to define exactly any unique lock-in frequency''}
\cite[p.70]{Gardner-1979-book}, \cite[p.188]{Gardner-2005-book}.

Below we demonstrate how to overcome these problems and rigorously define
a unique lock-in frequency and range.

We now consider a specific $\omega_{\Delta}^{\text{free}}$
and denote by $D_{\text{lock-in}}(\omega_{\Delta}^{\text{free}})$
the corresponding lock-in domain.
Such a domain exists for any
$|\omega_{\Delta}^{\text{free}}| \in \Omega_{\text{hold-in}}$
because at least the  equilibria are contained in this domain.
For a set $\omega_\Delta^{\text{free}} \in \Omega$
we consider the intersection of corresponding lock-in domains
(see, e.g. the intersections of local lock-in domains for various $\omega_{\Delta}^{\text{free}}=\pm|\widetilde{\omega}|$
in Fig.~\ref{lockin3def} ---  domains shaded both by red vertical and black horizontal lines):
\[
  {\rm D}_{\text{lock-in}}(\Omega) = \bigcap\limits_{\omega_\Delta^{\text{free}} \in \Omega} {\rm D}_{\text{lock-in}}(\omega_{\Delta}^{\text{free}}).
\]

\begin{definition} \label{def-lock}
 A lock-in range is the largest interval $[0,\omega_l)$
 such that for any $|\omega_\Delta^{\text{free}}| \in [0,\omega_l)$
 the mathematical model of the loop in the signal's phase space
 is globally asymptotically stable (i.e. $[0,\omega_l) \subset [0,\omega_p)$)
 and the following domain
 \[
  {\rm D}_{\text{lock-in}}\big((-\omega_l,\omega_l)\big) =
  \bigcap\limits_{|\omega_\Delta^{\text{free}}| < \omega_l} {\rm D}_{\text{lock-in}}(\omega_{\Delta}^{\text{free}}).
 \]
 contains all corresponding equilibria:
 \[
 \big( x_{eq}(\omega_{\Delta}^{\text{free}}),
     \theta_{eq}(\omega_{\Delta}^{\text{free}})\big)
 \in {\rm D}_{\text{lock-in}}\big((-\omega_l,\omega_l)\big).
 \]
\end{definition}
We call such domain ${\rm D}_{\text{lock-in}}={\rm D}_{\text{lock-in}}\big((-\omega_l,\omega_l)\big)$
\emph{a uniform lock-in domain} (uniform with respect to $(-\omega_l,\omega_l)$),
$\omega_{l}$ is called a \emph{lock-in frequency} (see \cite[p.40]{Gardner-1966}).
\smallskip

Various additional requirements may be imposed on the shape of the uniform lock-in domain
${\rm D}_{\text{lock-in}}$, e.g. it has to contain the line defined
by $x \equiv 0$ (see, e.g. \cite[p.258]{Kroupa-2012}) or the band defined by $|x| < c_{\text{max}}$.
If instead of global stability in the definition of the pull-in set
we consider stability in
the domain defined by $X_{\text{real}}$,
then we require that the intersection ${\rm D}_{\text{lock-in}} \bigcap X_{\text{real}}$
contains all corresponding equilibria.

\begin{remark}
  In the general case when there is no symmetry with respect to
  $\omega_{\Delta}^{\text{free}}$
  we have to consider a unsymmetrical interval containing zero
  in Definition~\ref{def-lock}.
\end{remark}

\noindent Similarly, we can define an extension of the lock-in range: $\Omega_{\text{lock-in}} \supset [0, \omega_l)$,
called a lock-in set (however, in general, such an extension may be not unique).

In other words, the definition implies that
\emph{if the loop is in a locked state,
then after an abrupt change of $\omega_{\Delta}^{\text{free}}$
within a lock-in range $[0, \omega_l)$,
the corresponding acquisition process in the loop leads,
if it is not interrupted, to a new locked state
without cycle slipping}.

Finally, our definitions give
\(
  \Omega_{\text{lock-in}} \subset \Omega_{\text{pull-in}}
  \subset \Omega_{\text{hold-in}}.
\)
If there is a certain stable equilibrium varies continuously
when $\omega_\Delta^{\text{free}}$ is changed within the
hold-in, pull-in, and lock-in ranges (see Footnote~\ref{fnote}), then
\[
  [0,\omega_l) \subset [0,\omega_p)
  \subset [0,\omega_h)
\]
which is in agreement with the classical consideration (see, e.g. \cite[p.34]{Best-1984},\cite[p.612]{HsiehH-1996},\cite[p.61]{Best-2007},\cite[p.138]{Egan-2007-book},\cite[p.258]{Kroupa-2012}).

\subsection{Approximations of the lock-in range of the classical PLL}
  For the case of the classical odd PD characteristic (see Fig.~\ref{lockin3def}),
  taking into account that equilibria are proportional
  to the frequency deviation (see \eqref{zeros1})
  and using the symmetry
  $\big(x_{eq}(\omega_{l}), \theta_{eq}(\omega_{l})\big) =
  - \big(x_{eq}(-\omega_{l}), \theta_{eq}(-\omega_{l})\big)$,
  we can effectively determine $\omega_{l}$.
  For that, we have to increase the frequency deviation
  $|\omega_{\Delta}^{\text{free}}|$  step by step
  and at each step, after the loop achieves a locked state,
  to change $\omega_{\Delta}^{\text{free}}=\widetilde{\omega}$ abruptly to $\omega_{\Delta}^{\text{free}}=-\widetilde{\omega}$
  and to check if the loop can achieve a new locked state without cycle
  slipping. If so, then the considered value $|\omega_{\Delta}^{\text{free}}|$
  belongs to $\Omega_{\text{lock-in}}$.
  If $\omega_{\Delta}^{\text{free}}\!=\!0$ belongs to $\Omega_{\text{pull-in}}$,
  then it is clear that $0$ belongs to
  $\Omega_{\text{lock-in}}$ (see Fig.~\ref{lockin3def},\,left).
  The limit value $\omega_{l}$ is defined by the case in Fig.~\ref{lockin3def},\,middle.
  At the next step
  when a value $|\omega_{\Delta}^{\text{free}}|=|\widetilde{\omega}| > \omega_{l}$ is considered,
  for $\omega_{\Delta}^{\text{free}}=-|\widetilde{\omega}|$
  the \newline
  trajectory from the initial point, corresponding to a stable equilibrium
  for $\omega_{\Delta}^{\text{free}}\!=\!|\widetilde{\omega}|$
  (see Fig.~\ref{lockin3def},\,right: red trajectory outgoing from a black dot),
  is attracted to an equilibrium
  only after cycle slipping.
  In other words \cite{KuznetsovLYY-2015-IFAC-Ranges}, for this case:

{\it
The lock-in range is a subset of the pull-in range
such that for each corresponding frequency deviation
the lock-in domain
(i.e. a domain of the loop states, where fast acquisition without cycle slipping is possible)
contains both symmetric locked states
(i.e. locked states for
\begin{figure}[ht]
 \centering
 \includegraphics[width=0.45\textwidth]{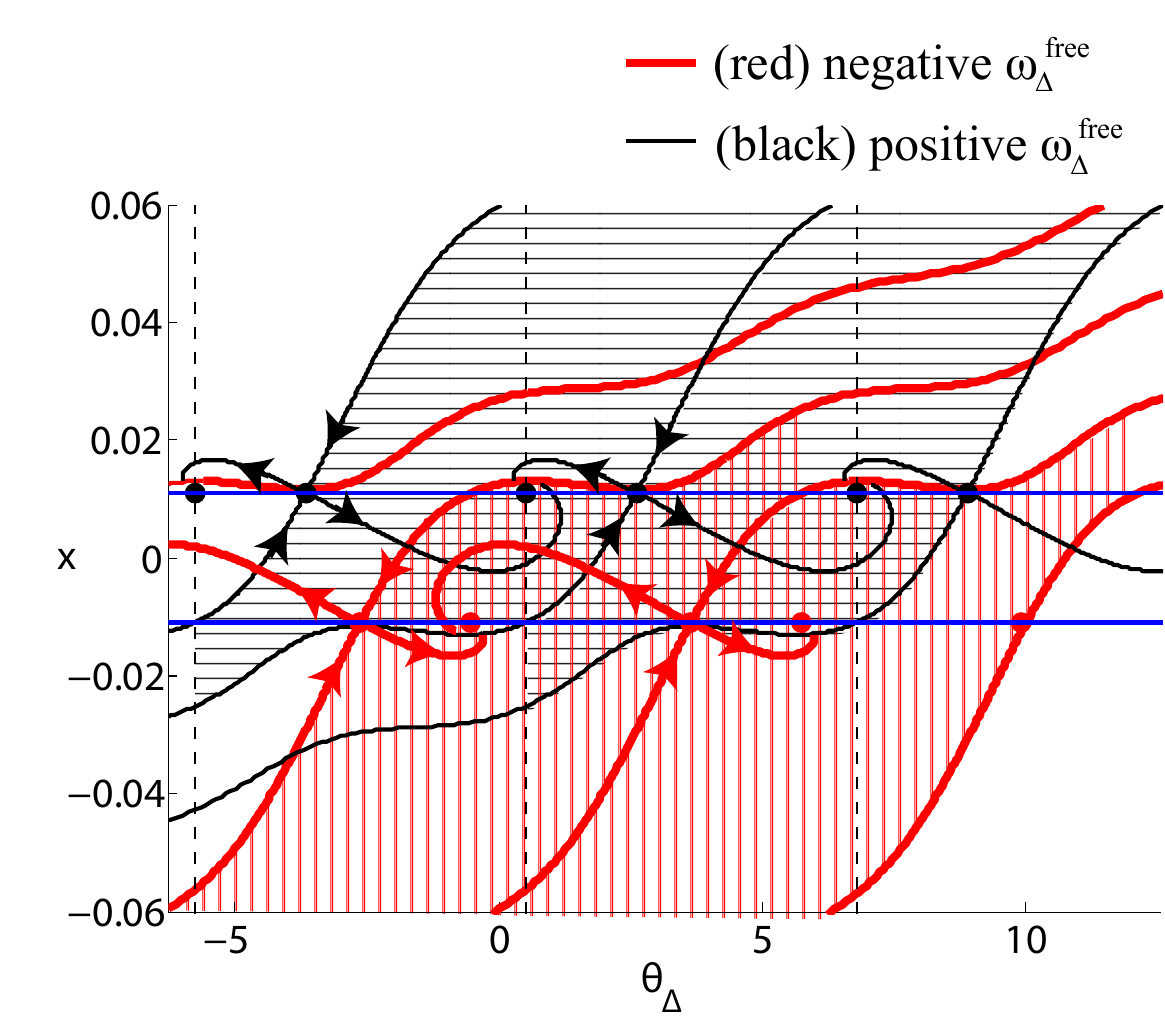}
 \caption{
 Phase portrait. Separatrices, equilibria
 and corresponding local lock-in domains (shaded):
 upper black is for $\omega_\Delta^{\text{free}} = 61.5$,
 lower red is for $\omega_\Delta^{\text{free}} = -61,5$.
 The uniform lock-in domain is approximated by the band between
 two blue horizontal lines: $|x|\leq 0.0110$.
 }. \label{lockindefband}
\end{figure}
the positive and negative value of the difference between the reference frequency and the VCO free-running frequency).
}

In Fig.~\ref{lockin3def},\,middle the set ${\rm D}_{\text{lock-in}}$:
contains all equilibria $x_{eq}(\omega_{\Delta}^{\text{free}})$ for
$0 \leq |\omega_{\Delta}^{\text{free}}| < \omega_{l}$.
 However for some non-equilibrium initial states
 from the band defined by $\{x: |x|< |x_{eq}(\omega_{l})|\}$
 (phase error $\theta_{\Delta}$ takes all possible values),
 cycle slipping can take place.
  For example, see the points to the left and to the right of the black equilibrium states
  (i.e. for $\omega_{\Delta}^{\text{free}}=|\omega_{l}|>0$),
  lying above the red separatrix (i.e. for $\omega_{\Delta}^{\text{free}}=-|\omega_{l}|<0$),
  correspond to the red trajectories
  (i.e. for $\omega_{\Delta}^{\text{free}}=-|\omega_{l}|<0$),
  which are attracted to an equilibrium only after cycle slipping.
  To approximate the ${\rm D}_{\text{lock-in}}$ by a band,
  $\omega_{l}$ can be slightly decreased to cut the above points.
  In Fig.~\ref{lockindefband}
  the band defined by $X_{\text{lock-in}}=\{x: |x|< |x_{eq}(\widetilde\omega_{l})|,\ \widetilde\omega_{l}< \omega_{l} \}$
  is contained in ${\rm D}_{\text{lock-in}}$
  and for any initial state from the band
  the corresponding acquisition process in the loop leads, if it is not interrupted,
  to lock up without cycle slipping.
  Such a construction is more laborious and
  requires rigorous analysis of the phase space or exhaustive simulation.

\begin{remark}
  If we define (see, e.g. \cite[p.92]{PurkayasthaS-2015})
  cycle slipping by the interval of maximum length $2\pi$
  instead of $4\pi$ in Definition~\ref{def-cs}:
  i.e. $\limsup_{t\to\infty} |\theta_\Delta(0) - \theta_\Delta(t)| > \pi$,
  then for any $|\omega_\Delta^{\text{free}}| > 0$
  a distance between neighboring unstable and stable equilibria
  and a phase deviation of the corresponding unstable saddle separatrix
  may exceed $\pi$ (see, e.g. Fig.~\ref{lockindefband}).
  Thus, the lock-in range may contain only $|\omega_\Delta^{\text{free}}| = 0$.
\end{remark}
\smallskip

\begin{remark}\label{remarkPI}
If the filter -- perfect integrator can be implemented in considered architecture,
the loop can be designed with the first order PI filter having
the transfer function $H(s) = \frac{1+s\tau_2}{s\tau_1}$.
Equations of the loop in this case become
\begin{figure}[ht]
 \centering
 \includegraphics[width=0.45\textwidth]{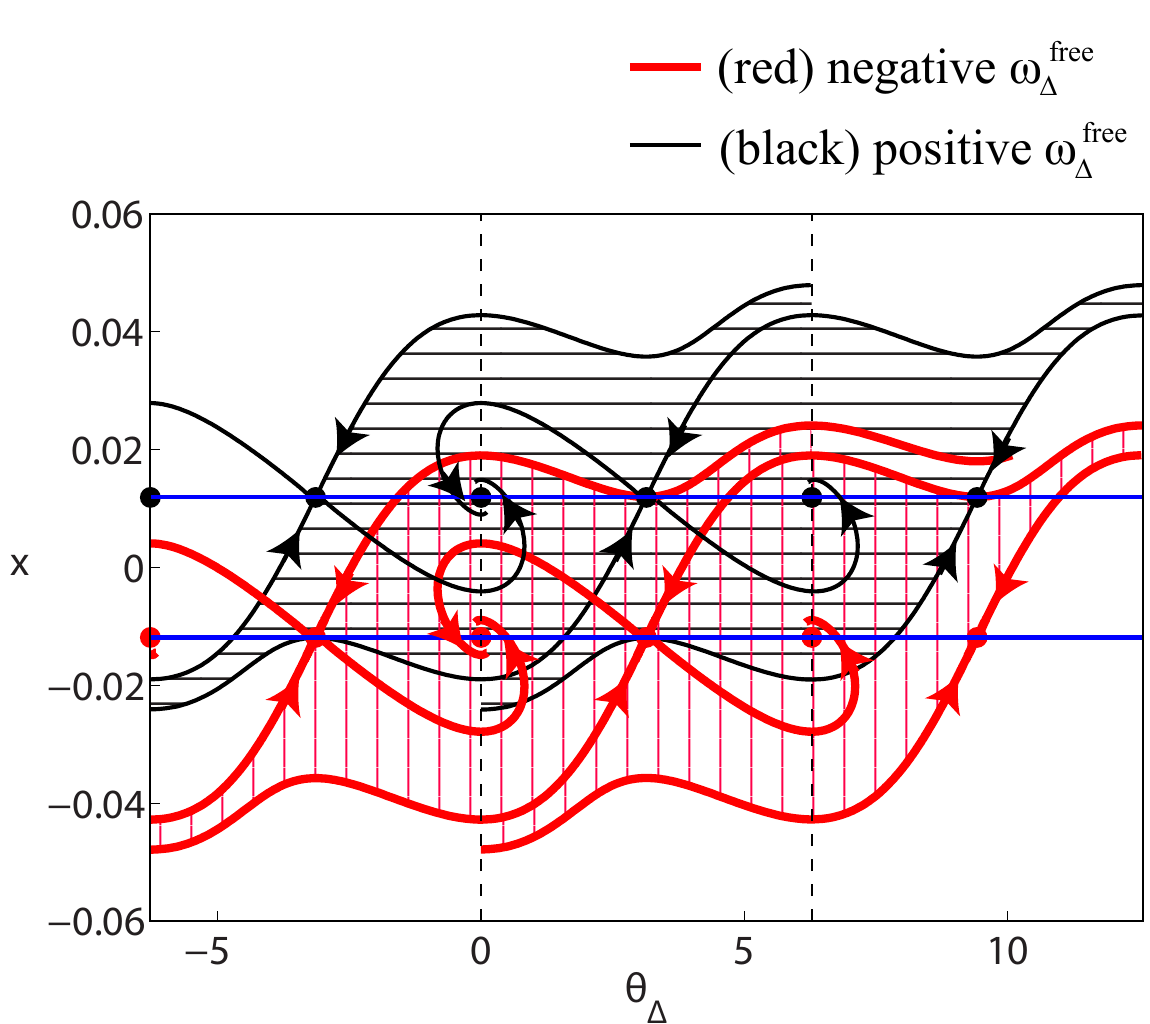}
 \caption{
 Phase portraits for the classical PLL with the following parameters:
 $H(s)=\frac{1+0.0225s}{0.0633s}$, $L=250$, and $\omega_\Delta^{\text{free}} = \pm 47$.
 Separatrices, equilibria
 and corresponding local lock-in domains (shaded):
 upper black is for $\omega_\Delta^{\text{free}} = 47$,
 lower red is for $\omega_\Delta^{\text{free}} = -47$.
 The uniform lock-in domain is approximated by the band between
 the two blue horizontal lines: $|x|\leq 0.0119$.
 } \label{lockindefband_pi}
\end{figure}
\begin{equation}\label{pi-system}
   \dot{x} = \frac{1}{\tau_1} \varphi(\theta_\Delta), \ \
   \dot\theta_\Delta = \omega_\Delta^{\text{free}}
   - Lx - L\frac{\tau_2}{\tau_1}\varphi(\theta_\Delta),
\end{equation}
or equivalently
\begin{equation}
\label{pi_equation}
  \ddot\theta_\Delta = -L\frac{1}{\tau_1} \varphi(\theta_\Delta)
  -L\frac{\tau_2}{\tau_1}\varphi'(\theta_\Delta)\dot\theta_\Delta.
\end{equation}
Here the equilibria are defined from the equations
\[
 \varphi(\theta_{eq}) = 0, \
 \ x_{eq} = {\omega_{\Delta}^{\text{free}} L^{-1}}.
\]
Because model \eqref{pi_equation} does not depend explicitly on $\omega_{\Delta}^{\text{free}}$,
the hold-in and pull-in ranges are either infinite or empty.
Note, that the parameter $\omega_{\Delta}^{\text{free}}$ shifts the phase plane vertically
(in the variable $x$) without distorting trajectories,
which simplifies the analysis of the uniform lock-in domain and range (see Fig.~\ref{lockindefband_pi}).
If the transfer function $H(s)$ of a high order filter
has the term $s^r$ with $r\in\mathbb{N}$ in the denominator,
then instead of equilibria we have a stationary linear manifold:
$
   \varphi(\theta_{eq}) = 0, \ c_1 x_{eq}^1 +\ldots+ c_r x_{eq}^r
   = \frac{-\omega_\Delta^{\text{free}}}{L}.
$

For the classical PLL with the filter's transfer function
$H(s) = \frac{\beta+\alpha s}{s}$
it can be analytically proved that the pull-in range is theoretically infinite.
Some needed explanations are given by Viterbi \cite{Viterbi-1966}
using phase plane analysis.
But, even in such a simple case, rigorous phase plane analysis
is a complex task
(e.g. \cite{AlexandrovKLNS-2015-IFAC},
the proof of the nonexistence of heteroclinic and first-order cycles is omitted in \cite{Viterbi-1966}).
The rigorous analytical proof can be effectively achieved
by considering a special Lyapunov function \cite{Bakaev-1963,LeonovK-2014,AlexandrovKLNS-2015-IFAC}: $V(x,\theta_\Delta)=\frac{1}{2}\big(x-\frac{\omega_\Delta^{\text{free}}}{L}\big)^2+\frac{2\beta}{L}\sin^2{\frac{\theta_\Delta}{2}} \geq 0$
and $\dot V(x,\theta_\Delta) = - h\beta\sin^2{\theta_\Delta} \leq 0$.
Here it is important that for any $\omega_\Delta^{\text{free}}$
the set $\dot V(x,\theta_\Delta)\equiv0$ does not contain the whole trajectories of system \eqref{pi-system} except for equilibria.
\end{remark}

\subsection{Initial and free-running frequencies of VCO}

Note that in the above Definitions~\ref{def-hold}, \ref{def-pull},
and \ref{def-lock} the hold-in, pull-in, and lock-in sets
are defined by the frequency deviation,
i.e. by the absolute value of the difference between
VCO free-running frequency (in the open loop)
and the input signal's frequency:
$|\omega_{\Delta}^{\text{free}}| = |\omega_1-\omega_2^{\text{free}}|$.
The VCO free-running frequency $\omega_2^{\text{free}}$
is different from the VCO initial frequency $\omega_2(0)$:
\(   \omega_2(0) = \omega_2^{\text{free}} + g(0), \)
 where $g(0) = c^{*}x(0)+h\varphi(\theta_{\Delta}(0))$
 is the initial control signal,
 depending on the initial states of the filter $x(0)$
 and the initial phase difference $\theta_{\Delta}(0)$.

 It is interesting that for simplified model \eqref{mathmodel-class-simple}
 with $h=0$ (see eq.~2.20 in the classic reference \cite{Viterbi-1966})
 the absolute value of the initial difference between frequencies
 $|\dot \theta_{\Delta}(0)| = |\omega_{\Delta}(0)|=|\omega_1-\omega_2(0)|$ is equal
 to the frequency deviation $|\omega_{\Delta}^{\text{free}}| = |\omega_1-\omega_2^{\text{free}}|$.
 Following such simplified consideration
 in engineering literature the concept of an ``{\it initial frequency difference}''
 can be found to be in use instead of the concept  of a ``{\it frequency deviation}'':
 see, e.g. \cite[p.44]{Gardner-1966}
 ``{\it If the initial frequency difference (between VCO and input)
 is within the pull-in range, the VCO frequency will slowly change
 in a direction to reduce the difference}'',
  \cite[p.1792]{Chen-2002-book} ``{\it The maximum frequency difference
  between the input and the output that the PLL can lock
 within one single beat note is called the lock-in range of the PLL}'',
  \cite[p.49]{KiharaOE-2002} ``{\it  Whether the PLL can get synchronized
 at all or not depends on the initial frequency difference between
 the input signal and the output of the controlled oscillator.}''
 In general, the change of $\omega_2^{\text{free}}$ to $\omega_2(0)$
 may lead to wrong results in the above definitions of ranges
 because in the case of $x(0) \neq 0$, $h \neq 0$
 or non-odd function $\varphi(\theta_{\Delta})$
 for the same values of $\omega_2(0)$
 the loop can achieve synchronization or not
 depending on the filter's initial state $x(0)$,
 the initial phase difference $\theta_{\Delta}(0)$, and $\omega_2^{\text{free}}$.
 See the corresponding example.

\begin{example}\label{change}
 Consider the behavior of model \eqref{final_system}
 for the sinusoidal signals
 (i.e. $\varphi(\theta_{\Delta}) = {1 \over 2}\sin(2\theta_{\Delta})$)
 and the fixed parameters:
 $\omega_\Delta = 100, H(s) = \frac{(1+s \tau_2)}{1+s(\tau_1 + \tau_2)},
 \tau_1 = 0.0448, \tau_2 = 0.0185,
 L = 250$.
 In Fig.~\ref{small omega vatiation sim3}
 the phase portrait of system \eqref{final_system} is shown.
 The blue dash line consists of points for which the initial frequency difference
 is zero: $\omega_\Delta(0)=\dot\theta_\Delta(0)=0$.
 Despite the fact that the initial frequency differences of all trajectories outgoing
 from the blue line are the same (equal to $0$),
 the green trajectory tends to a locked state
 while the magenta trajectory can not achieve this.
\begin{figure}[h]
 \centering
 \includegraphics[width=0.45\textwidth]{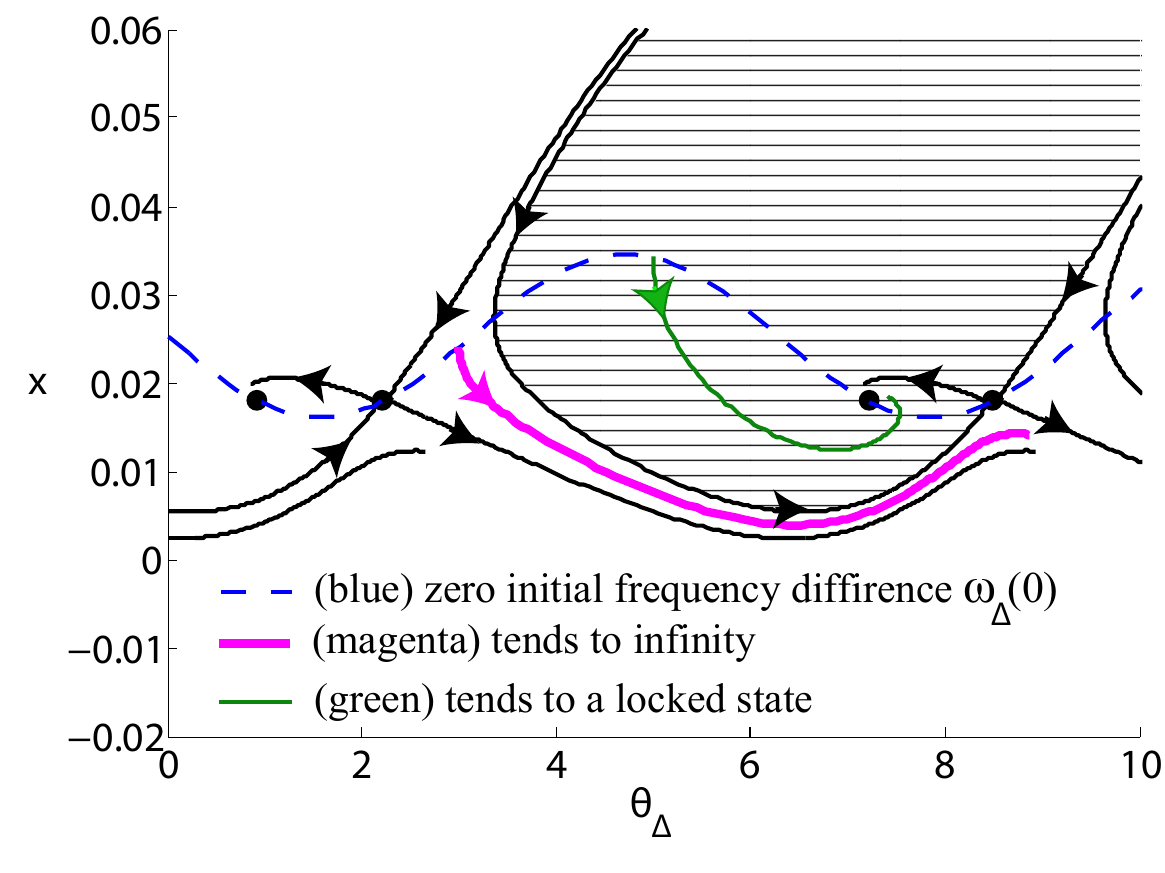}
 \caption{
  Phase portrait for $\omega_{\Delta}^{\text{free}} = 100$.
  Blue dash curve corresponds to the set defined by $\dot\theta_\Delta(0)=0$.
  Initial points of the green (upper) and magenta (lower) trajectories correspond
  to the same initial frequency difference $\omega_\Delta(0)=0$.}
 \label{small omega vatiation sim3}
 \end{figure}
 \end{example}

\section{\uppercase{Conclusions}}
This survey discussed a disorder and inconsistency in the definitions of ranges
currently used.
An attempt is made to discuss and fill some of the gaps identified
between mathematical control theory,  the theory of dynamical systems
and the engineering practice of phase-locked loops.
Rigorous mathematical definitions for
the hold-in, pull-in, and lock-in ranges are suggested.
The problem of unique lock-in frequency definition,
posed by Gardner \cite{Gardner-1979-book},
is solved and an effective way to determine
the unique lock-in frequency is suggested.

\section*{\uppercase{Acknowledgements}}
 This work was supported by the Russian Scientific Foundation (project 14-21-00041)
 and Saint-Petersburg State University. %6.39.416.2014 (sections I-IV)
 The authors would like to thank Roland~E.~Best,
 the founder of the Best Engineering Company, Oberwil, Switzerland
 and the author of the bestseller on PLL-based circuits \cite{Best-2007}
 for valuable discussion. % on PLL-based circuits.

%\bibliographystyle{apalike}
%\bibliographystyle{IEEEtran}
%\bibliography{../../bib/bib_pll,../../bib/bib_nk,../../bib/bib_leonov,../../bib/bib_full,../../bib/bib-gly}
%\bibliography{C:/Dropbox/bib/bib_nk,C:/Dropbox/bib/bib_leonov,C:/Dropbox/bib/bib_full,C:/Dropbox/bib/bib_pll,C:/Dropbox/bib/bib-gly}

\newcommand{\noopsort}[1]{} \newcommand{\printfirst}[2]{#1}
  \newcommand{\singleletter}[1]{#1} \newcommand{\switchargs}[2]{#2#1}

%\vspace{-1.5cm}
\begin{IEEEbiography}
 [{\includegraphics[width=1in,height=1.25in,clip,keepaspectratio]{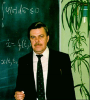}}]
{Gennady Leonov}
received his Candidate degree in 1971 and Dr.Sci. in 1983 from Saint-Petersburg State University.
In 1986 he was awarded the USSR State Prize
\emph{for development of the theory of phase synchronization for radiotechnics and communications}.
Since  1988 he has been Dean of the Faculty of Mathematics and Mechanics
at Saint-Petersburg State University
and since 2007 Head of the Department of Applied Cybernetics.
%Professor G.A. Leonov was awarded the Prize of St.-Petersburg State University (1985),
%State Prize of USSR (1986), the Prize of Technische Universitet Dresden (1990),
%Medal of the University of Jyv\"{a}skyl\"{a} (2011),
%Andronov prize of the Russian Academy of Science (2012).
He is member (corresponding) of the Russian Academy of Science, in 2011 he was elected to the IFAC Council.
His research interests are now in control theory and dynamical systems.
\end{IEEEbiography}
%\vspace{-1.6cm}

\begin{IEEEbiography}
 [\vspace{-0.8cm}{\includegraphics[width=1in,height=1.25in,clip,keepaspectratio]{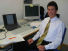}}]
 {Nikolay Kuznetsov}
 received his Candidate degree from Saint-Petersburg State University
 (2004) and PhD from the University of Jyv\"{a}skyl\"{a} (2008).
 He is currently Deputy Head of the Department of Applied Cybernetics
 at Saint-Petersburg State University and
 Adjunct Professor at the University of Jyv\"{a}skyl\"{a}.
 His interests are now in dynamical systems stability and oscillations,
 Lyapunov exponent, chaos, hidden attractors, phase-locked loop nonlinear analysis,
 nonlinear control systems.
  \newline
 E-mail: nkuznetsov239@gmail.com (corresponding author)
\end{IEEEbiography}
%\vspace{-1.4cm}
\begin{IEEEbiography}
 [\vspace{-0.8cm}{\includegraphics[width=1in,height=1.25in,clip,keepaspectratio]{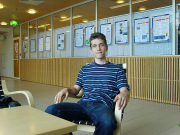}}]
{Marat Yuldashev}
received his Candidate degree from St.Petersburg State University
(2013) and PhD from the University of Jyv\"{a}skyl\"{a} (2013).
He is currently at Saint-Petersburg University.
His research interests cover nonlinear models of phase-locked loops and Costas loops, and SPICE simulation.
\end{IEEEbiography}
%\vspace{-1.4cm}
\begin{IEEEbiography}
 [\vspace{-0.8cm}{\includegraphics[width=1in,height=1.25in,clip,keepaspectratio]{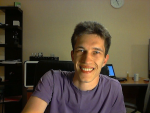}}]
{Renat Yuldashev}
received his Candidate degree from St.Petersburg State University
(2013) and PhD from the University of Jyv\"{a}skyl\"{a} (2013).
He is currently at Saint-Petersburg University.
His research interests cover nonlinear models of phase-locked loops and Costas loops,
and simulation in MatLab Simulink.
\end{IEEEbiography}
\end{document}